\documentclass[review,onefignum,onetabnum]{siamart171218}

\usepackage{lipsum}
\usepackage{amsfonts}
\usepackage{graphicx}
\usepackage{epstopdf}
\usepackage{algorithmic}
\usepackage{hyperref}
\usepackage{amssymb}
\usepackage{amsmath}
\usepackage{float}
\usepackage{color}
\graphicspath{{./Images/}}

\ifpdf
  \DeclareGraphicsExtensions{.eps,.pdf,.png,.jpg}
\else
  \DeclareGraphicsExtensions{.eps}
\fi

\newsiamremark{remark}{Remark}
\newsiamremark{hypothesis}{Hypothesis}
\crefname{hypothesis}{Hypothesis}{Hypotheses}
\newsiamthm{claim}{Claim}
\newsiamthm{exmp}{Example}

\headers{FSI for the Classroom: Convergence!}{N.A. Battista \& M. S. Mizuhara}

\title{Fluid-Structure Interaction for the Classroom: Speed, Accuracy, Convergence, and Jellyfish! \thanks{Submitted to the editors \textcolor{red}{February 20, 2019}.
\funding{This work was funded by NSF OAC-1828163 and the Support of Scholarly Activities Grant from TCNJ (The College of New Jersey)}}}

\author{Nicholas A. Battista\thanks{Department of Mathematics and Statistics, The College of New Jersey, 2000 Pennington Road, Ewing Township, NJ 08628; corresponding author
  (\email{battistn@tcnj.edu}, \url{http://battistn.pages.tcnj.edu/}).}
\and
Matthew S. Mizuhara\thanks{Department of Mathematics and Statistics, The College of New Jersey, 2000 Pennington Road, Ewing Township, NJ 08628 
  (\email{mizuharm@tcnj.edu}, \url{http://www.tcnj.edu/\%7Emizuharm/}).}}

\usepackage{amsopn}

\ifpdf
\hypersetup{
  pdftitle={Fluid-Structure Interaction for the Classroom: Speed, Accuracy, Convergence, and Jellyfish!},
  pdfauthor={N.A. Battista and M.S. Mizuhara}
}
\fi

\begin{document}

\maketitle

%%%%%%%%%%%%%%%%%%%%%%%%%%%%%%%%%%%%%%%%%%%%%%%%%%%%%%%%%%%%%%%%%%%%%%%
%
% ABSTRACT
%
%%%%%%%%%%%%%%%%%%%%%%%%%%%%%%%%%%%%%%%%%%%%%%%%%%%%%%%%%%%%%%%%%%%%%%%

\begin{abstract}

\textit{When is good, good enough?} This question lingers in approximation theory and numerical methods as a competition between accuracy and practicality. Numerical Analysis is traditionally where the rubber meets the road: students begin to use numerical algorithms to compute approximate solutions to non-trivial problems. However, it is difficult for students to understand that more accuracy is not always the goal, but rather \textit{enough} accuracy for practical use and meaningful interpretation. %Moreover, it is difficult for students to reconcile that the practical goal is not mere accuracy but rather \textit{enough} accuracy for meaningful interpretation.
This compromise between accuracy and computational time/resources can be explored through the use of convergence plots. We offer a variety of classroom activities that allow students to discover the usefulness of convergence plots, including a contemporary example involving jellyfish locomotion using fluid-structure interaction modeling. These examples will additionally illustrate subtleties of convergence analysis: the same numerical scheme can exhibit different convergence rates, and the definition of ``error'' may change the convergence properties. To solve the fluid-structure interaction system, the open source software \textit{IB2d} is used.  All numerical codes, scripts, and movies are provided for streamlined integration into a classroom setting.

\end{abstract}

% REQUIRED

% REQUIRED

% 65-01: numerical analysis exposition
% 65D05: interpolation
% 65D07: Splines
% 97M10: Modeling and interdisciplinarity
% 97M60: Biology, chemistry, medicine
% 97N10: Comprehensive works
% 97N40: Numerical analysis (math ed)
% 97N50: Interpolation and approximation (math ed)
% 97N80: Mathematical software, computer programs (math ed)
% 76M25: Fluid mechanics Other numerical methods
% 76Z10: Biopropulsion in water and in air
% 76Z99: Biofluids, None of the above, but in this section
% 92C10: Biomechanics

%%%%%%%%%%%%%%%%%%%%%%%%%%%%%%%%%%%%%%%%%%%%%%%%%%%%%%%%%%%%%%%%%%%%%%%
%
% INTRODUCTION
%
%%%%%%%%%%%%%%%%%%%%%%%%%%%%%%%%%%%%%%%%%%%%%%%%%%%%%%%%%%%%%%%%%%%%%%%

%%%%%%%%%%%%%%%%%%%%%%%%%%%%%%%%%%%%%%%%%%%%%%%%
%
%
% INTRODUCTION
%
%
%%%%%%%%%%%%%%%%%%%%%%%%%%%%%%%%%%%%%%%%%%%%%%%%

\section{Introduction}
\label{intro}

There is often an internal struggle for mathematics students to turn to computers to \textit{approximate} solutions to non-trivial problems. It seems to them a challenge to the belief that the beauty of mathematics comes from its precision, rigidity, and consistency. Opening themselves to a new tool that openly boasts ``I can live with this amount of error" seems to contradict a lot of students' maturing mathematical identities, or in the very least cause an unsettled feeling. 

The question of ``when is good, good enough?" also poses an interesting philosophical question for students. Going from a world of precise solutions to one of approximations usually makes many students immediately jump towards approximations with the highest accuracy possible. This is not a bad endeavor; new numerical techniques are constantly being developed and discovered that push towards higher accuracy. However, the subtlety is that for \textit{practical} purposes these techniques must also boast that they don't require immense additional computational effort, that is, they do not take longer than previously implemented methods or require specific computational infrastructure or resources that may not be available. 

\begin{figure}[H]
    %\centering
    \centering
    \includegraphics[width=0.99\textwidth]{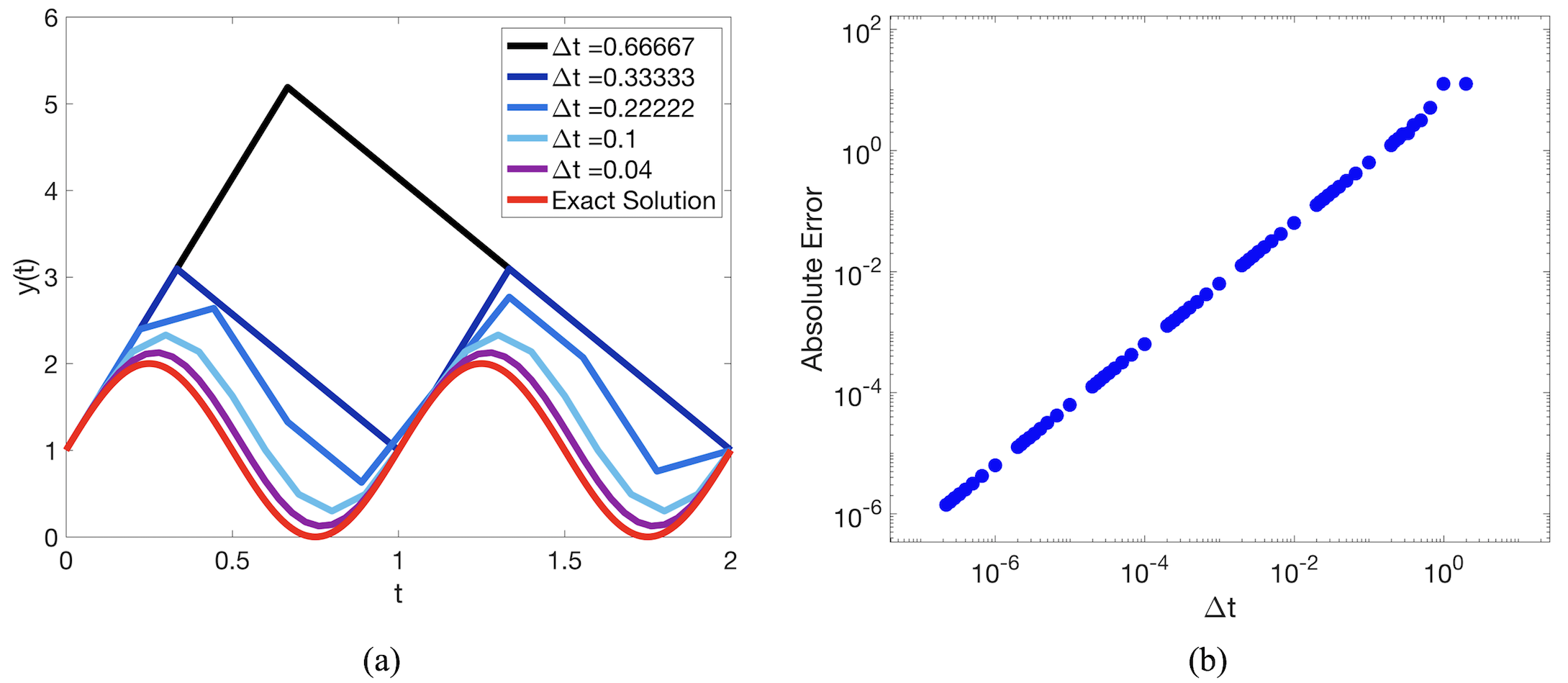}
    \caption{(a) Showing increasing accuracy of numerical solutions vs. the exact solution of \eqref{eq:euler1}-\eqref{eq:euler2} by increasing the time-step resolution, i.e., decreasing $\Delta t$. (b) Convergence plot illustrating the absolute error between the numerical solution and exact solution for different $\Delta t.$ See Appendix \ref{app:euler} for a detailed discussion of the problem.}
    \label{fig:EulerComp}
\end{figure}

This competition between speed and accuracy can be demonstrated through the use of convergence plots. Figure \ref{fig:EulerComp} illustrates how for a particular problem, the numerical solutions achieve greater accuracy for smaller $\Delta t$. Figure \ref{fig:EulerComp}a shows numerical approximations of
\begin{align}
    \frac{dy}{dt} &= 2\pi \cos(2\pi t)\label{eq:euler1} \\
    y(0) &= 1\label{eq:euler2}
\end{align}
on $[0,2]$ using the forward Euler method and their qualitative closeness to the exact solution $y=\cos(2\pi t)$.
%qualitatively shows how close the numerical approximations are to the exact solution, while 
Figure \ref{fig:EulerComp}b gives the \textit{convergence plot}, that is, a study of how the absolute error decreases for different values of $\Delta t.$ The absolute error is defined to be the maximum of the absolute value of the difference between the numerical approximation and the exact solution, e.g., the $L^\infty$-distance. The absolute error vs. $\Delta t$ is given on a log-log plot, providing us information on the \textit{convergence rate} of the numerical scheme. It is evident that there is a linear relationship between the $\log(\mbox{absolute error})$ and $\log(\Delta t)$; we compute the slope of this line, call it $m$. Hence we have $$\log(\mbox{absolute error}) \sim m\log(\Delta t).$$

%\textcolor{red}{IS IT POSSIBLE THAT $\log(error)=m\log(\Delta t)+C$? OR SHOULD $C=0$ ALWAYS?}
%\textcolor{blue}{ Great point, Matt. To be rigorous, these should really all be asymptotic relations and not equations :)}
Solving this equation for the absolute error allows us to see how fast we expect the error to decay as a function of $\Delta t$, e.g.,

\begin{equation}
    \label{eq:eulerError} \mbox{absolute error} \sim \Delta t^m.
\end{equation}

We would then see that this numerical scheme applied to this particular problem is approximately $m^{th}$ order accurate. Note that $m$ will generally not be an integer. While one may be inclined  to round it to the nearest integer to get a better representation of an integer \textit{order of convergence}, some methods are historically known to have non-integer order of convergence, such as the secant method for approximating roots on a nonlinear equation, which formally has an order of convergence equal to the Golden Ratio, $\frac{1+\sqrt{5}}{2}\approx 1.618$, see Appendix \ref{app:secant_convergence}.

Students may be familiar with these plots to show how accurate a numerical scheme is, i.e., is it $1^{st}$, $2^{nd}$, or generally $m^{th}$ order, but behind these plots there is something else that is usually particularly subtle (and menacing) - the accuracy one can hope to achieve for a certain amount of computational time. Recall that convergence plots typically show some metric of error vs. resolution%, whether number of iterations performed or temporal or spatial resolutions
; however, they provide no direct information on how long it took a simulation to run. The story of required computational time for an algorithm is hidden in its resolution. That is, the computational time scales with the resolution! 

\begin{figure}[H]
    %\centering
    \centering
    \includegraphics[width=0.65\textwidth]{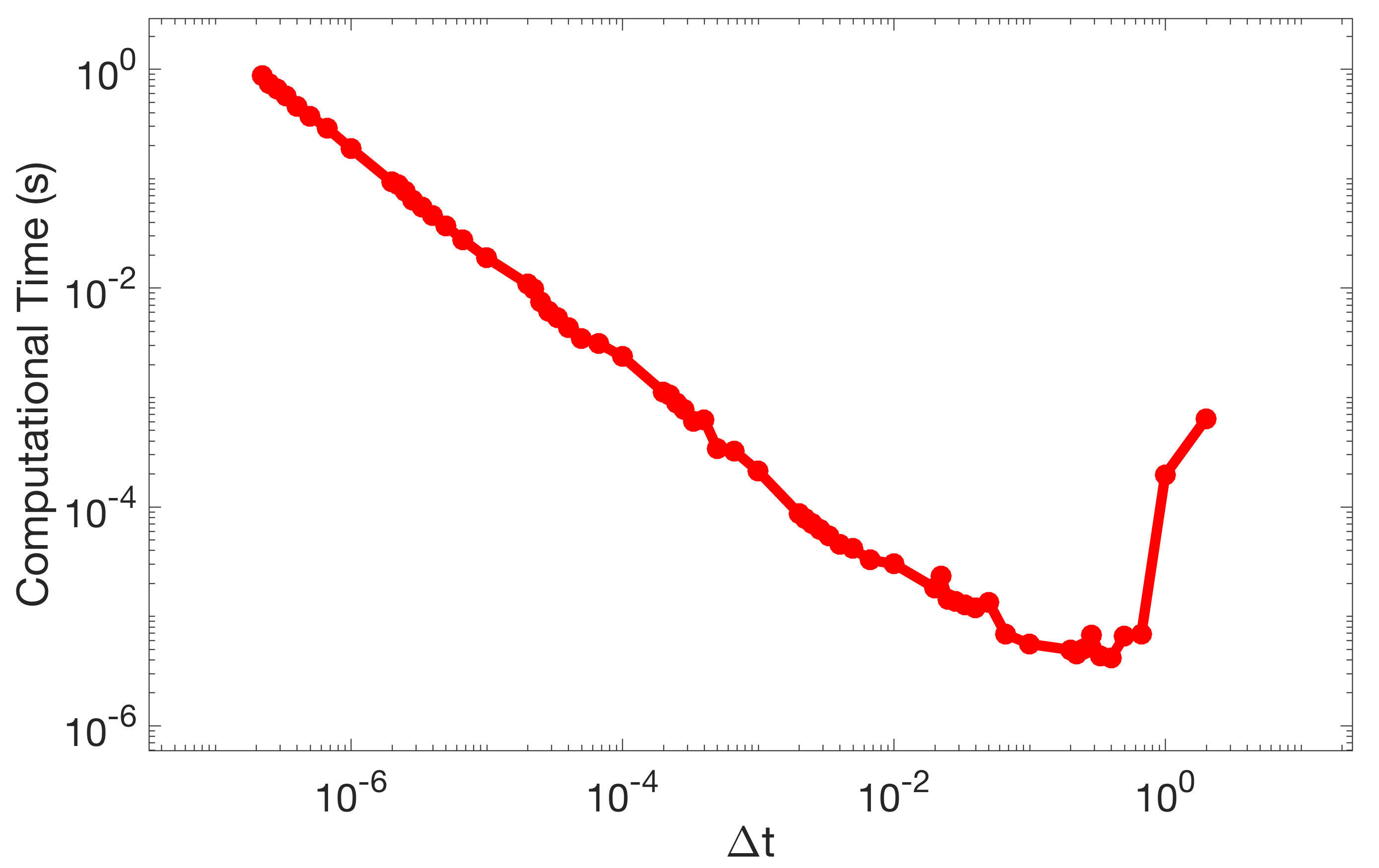}
    \caption{Provides the computational time vs. $\Delta t$ for simulation of \eqref{eq:euler1}-\eqref{eq:euler2}. As $\Delta t$ decreases (and the accuracy of the numerical approximation increases) the algorithm takes longer to run.}
    \label{fig:EulerTime}
\end{figure}

For example, Figure \ref{fig:EulerTime}, which gives the computational time vs. $\Delta t$ to numerically solve \eqref{eq:euler1}-\eqref{eq:euler2}, complements the convergence plot previously shown in Figure \ref{fig:EulerComp}. As $\Delta t$ decreases, the accuracy increases; however, this plot illustrates that the computational time required to run the algorithm also increases. Using the same rigmarole as when computing the order of convergence, we see that for this scheme that the computational time required also scales as $$\mbox{computational time} \sim \Delta t^{\mbox{slope}},$$
where in this case the slope is \textit{negative}. In this particular example, the computational time required increases inversely proportional with the time-step, $\Delta t$, e.g, computational time increases as time resolution increases. For a detailed discussion of the algorithm and problem in this example see Appendix \ref{app:euler}.

If the resolution of a numerical scheme is significantly increased to offer much higher accuracy (like reducing $\Delta t$ above), the time a simulation requires to run may also dramatically skyrocket as a result. What good is a numerical approximation that gives 14 decimal accuracy if it takes $1$ year for the simulation to run, while an approximation with 10 decimal accuracy takes 4 days, or with 8 decimals takes 1 hour? \footnote{Of course, this is just a generalized analogy and particular applications may warrant and demand certain accuracies, so we're constrained. This is also where the community of numerical analysts focus much of their efforts - developing faster and more accurate methods.}

In essence convergence plots allow us to designate some type of self-efficiency metric on our numerical method. Rather than seeking an answer to the question, ``how can I achieve the most accuracy?" it is appropriate to morph that question into, ``for a specific desired accuracy, what resolution can I use to find a solution in a \textit{timely and practical} manner?" Of course in different application contexts \textit{timely} could mean very different things - from fractions of a second to days to even a week, a month, or longer!

In this paper we offer multiple projects to acquaint students with the idea of a convergence plot - one on approximating the Golden Ratio, another on the accuracy and convergence of the composite trapezoid rule, and one that arises out of popular aquatic locomotion studies of jellyfish swimming. The jellyfish example comes from contemporary fluid-structure interaction research \cite{Hoover:2015}. From all exercises students will see the beauty, practicality, and importance of convergence plots. Moreover the jellyfish example offers students the opportunity to see the practicality of convergence plots at the frontier of research; it allows a unique class activity (or course project) that bridges the interface of modern computational science,  and mathematical modeling. 

For details regarding the fluid-structure interaction software, see Appendix \ref{IB:Appendix}, or \cite{Battista:2015,BattistaIB2d:2016,BattistaIB2d:2017} for a more detailed overview. All simulations presented here are available on \url{https://github.com/nickabattista/ib2d} and can found in the sub-directory \textit{IB2d/matIB2d/Examples/Examples$\_$Education/Convergence} as well as the Supplementary Materials.

%%%%%%%%%%%%%%%%%%%%%%%%%%%%%%%%%%%%%%%%%%%%%%%%%%%%%%%%%%%%%%%%%%%%%%%
%
% GOLDEN RATIO
%
%%%%%%%%%%%%%%%%%%%%%%%%%%%%%%%%%%%%%%%%%%%%%%%%%%%%%%%%%%%%%%%%%%%%%%%

%%%%%%%%%%%%%%%%%%%%%%%%%%%%%%%%%%%%%%%%%%%%%%%%
%
%
% GOLDEN RATIO
%
%
%%%%%%%%%%%%%%%%%%%%%%%%%%%%%%%%%%%%%%%%%%%%%%%%

\section{Convergence to the Golden Ratio}
\label{sec:golden_ratio}

The infamous Golden Ratio, $\phi$, has popped up into many unsuspecting places in nature, from seed heads, human faces, hurricane and galaxy formations, music, and the building blocks of life - DNA \cite{Tung:2007,Murali:2012}.  The story begins with the Fibonacci Sequence, $F_n$, defined recursively by 

\begin{equation}
    \label{eq:fibonnacci} F_{n+1} = F_{n} + F_{n-1}, \mbox{ for $n\geq 1$},
\end{equation}

 with $F_0= F_1 =1$. While the Golden Ratio has many derivations, we will define it to the the ratio of successive terms of the following sequence

\begin{equation}
    \label{eq:phi} \phi = \lim_{n\rightarrow\infty} \phi_n =  \lim_{n\rightarrow\infty}  \frac{F_{n+1}}{F_n}.
\end{equation}

%Dusting off second semester Calculus (or introductory real analysis), one can ask students if they remember criteria for convergence of a sequence. 
%Hopefully they recall when the sequence is monotonic and bounded. 
%We can prove that the the Golden Ratio Sequence, Eq.(\ref{eq:phi}), is bounded, see Appendix \ref{app:fibonacci}. Also 
Interestingly, the even terms in the Golden Ratio sequence form a monotonically increasing sequence, while the odd terms form a monotonically decreasing sequence \cite{Herrman:2012}. Although some explore the implications of these subsequences, we suffice to remark that the sequence is Cauchy and thus convergent, see Appendix \ref{app:fibonacci}.

Once it's known that the Golden Ratio sequence converges, one can compute its limit rather elegantly: starting with the recursive definition of the Fibonacci Sequence and dividing both sides by $F_n$ we see that,

$$\frac{F_{n+1}}{F_n} = 1 + \frac{F_{n-1}}{F_n}$$

and hence

$$\phi_{n} = 1 + \frac{1}{\phi_{n-1}}.$$

In the limit as $n\rightarrow\infty$, $\phi_{n-1},\phi_n \to \phi$, so the above expression becomes

\begin{equation}
    \label{eq:GR_Quadratic} \phi = 1 + \frac{1}{\phi} \Rightarrow \phi^2 - \phi - 1 = 0.
\end{equation}

Solving this quadratic gives $$\phi = \frac{1\pm\sqrt{5}}{2}.$$

Since all the terms of the sequence are positive, we take the positive root and find the Golden Ratio to be $$\phi = \frac{1+\sqrt{5}}{2} \approx 1.61803398874989485.$$

We could have instead calculated successive approximations to $\phi$ by simply computing terms of the Fibonacci Sequence and taking the ratio of successive terms appropriately. Why would anyone do this? Well, posed as a more ``numerical analysis-y'' question: \textit{how many terms do we need to go out in the sequence $\phi_n$ before our approximation to $\phi$ is accurate to $4$ decimal places? $8$ decimal places? $15?$}

Hopefully it is safe to say that this is where a computer comes in rather handy, instead of pen and paper. Also this will naturally lead us to the idea of a convergence plot. Luckily for us we have an ace up our sleeve since we already know the true value of $\phi$. To a standard $64$-bit computer we see its value actually takes approximately $\hat{\phi}=1.618033988749895$ due to finite-precision restrictions \cite{Kincaid:2002,Burden:2014}. 

Using this value of $\hat{\phi}$ as our ``exact'' value, we can ask how much error there is associated with a particular approximation to the Golden Ratio, $\phi_n$. We can compute the \textit{absolute error} between $\hat{\phi}$ and $\phi_n$, where the absolute error is given by 

\begin{equation}
    \label{eq:abs_err_GR} E_n = \left| \hat{\phi} - \phi_n \right|.
\end{equation}

\begin{figure}[H]
    %\centering
    \centering
    \includegraphics[width=0.65\textwidth]{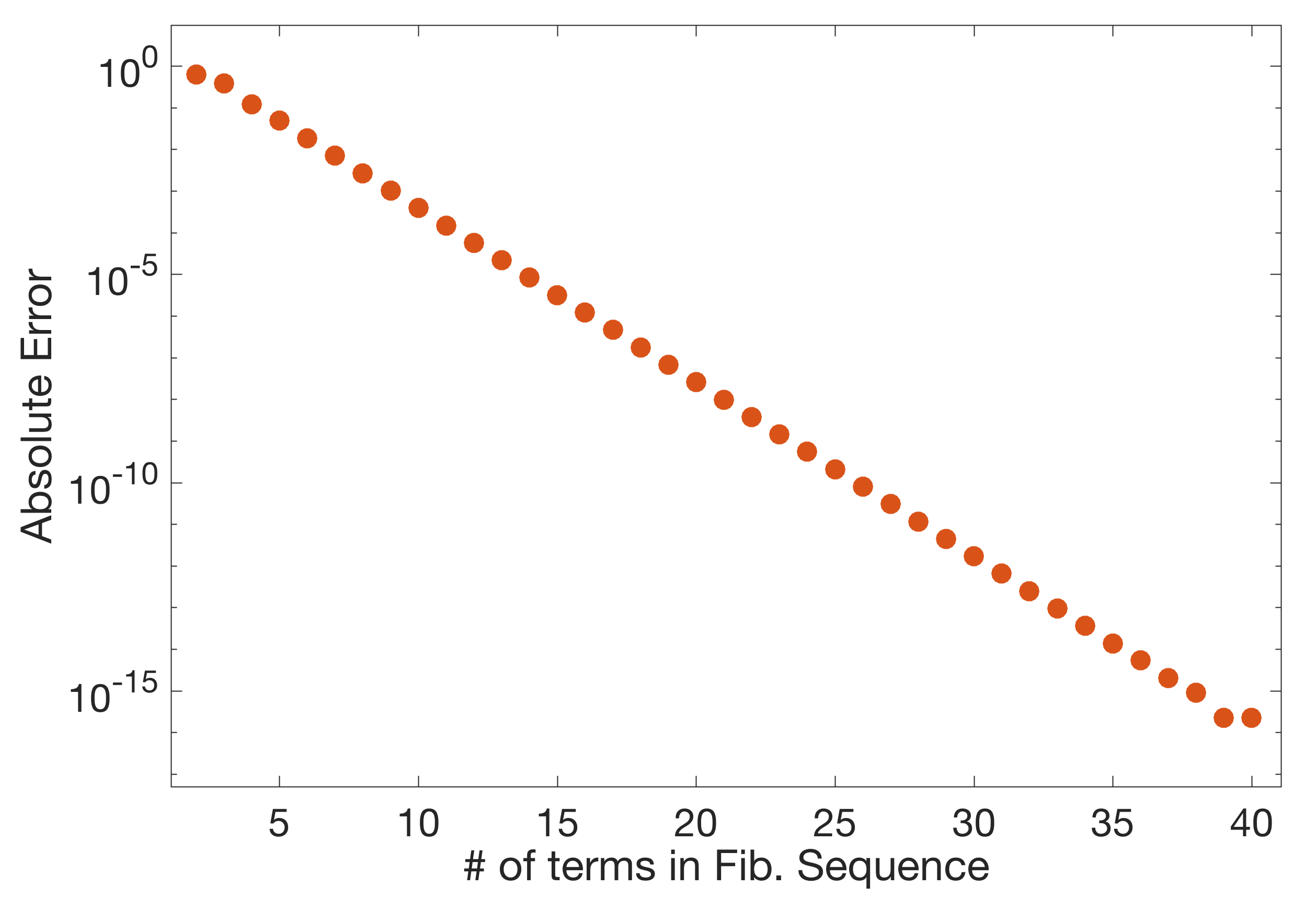}
    \caption{The absolute error for each approximation $\phi_n$, the $n^{th}$ approximation to the Golden Ratio, vs $n$, the number of the terms in the Fibonacci Sequence used for that approximation.}
    \label{fig:GoldenRatio_Conv}
\end{figure}

%Calculating each successive Golden Ratio approximation by including more terms in the Fibonacci Sequence and checking the accuracy of each approximation %leads to rather interesting information. It provides an algorithm to determine the accuracy of the $n^{th}$ Golden Ratio approximation, $\phi_n$. \textcolor{blue}{A script to do just this is found in\ldots}. 
 Accuracy for the first $40$ approximations is provided in Figure \ref{fig:GoldenRatio_Conv}. 

%\textcolor{red}{Out of curiosity, what is the slope of this line? My hunch is that it should converge something like $err = 10^{-\phi N}$. or maybe $err = 10^{-(\phi+\psi)N}$ where $\psi$ is the other root of $\phi^2-\phi-1=0$.}

Figure \ref{fig:GoldenRatio_Conv} is a \textit{convergence plot} showing the amount of error of as a function of the number of terms in the Fibonacci sequence. As more terms are including in the sequence, the Golden Ratio approximation becomes increasingly accurate. In fact, it looks as though the error decreases \textit{geometrically}, since the relationship between the logarithm of absolute error and number of terms appears linear, e.g., 
$$\log(\mbox{Absolute. Error}) \sim (\mbox{slope})N \ \ \Rightarrow \ \ \mbox{Absolute Error} \sim 10^{(\mbox{slope})N},$$
where $N$ is the number of terms in the Fibonacci Sequence and the slope is negative. Hence the error in successive Golden Ratio approximations decreases rapidly! This is called \textit{geometric} (or \textit{exponential}) convergence. Justification of this convergence is provided in Appendix \ref{app:fibonacci}.

Approximating the Golden Ratio with the first $40$ terms of the Fibonacci Sequence results in an error of approximately $10^{-15}$. However, it would not make a difference if we included more terms in the Fibonacci Sequence, as a standard 64-bit computer would not recognize any further digits. This is so-called \textit{machine precision}; we cannot get any more decimals of accuracy, without handing the floating point numbers in a special way \cite{Akkas:2003}. 

\begin{figure}[H]
    %\centering
    \centering
    \includegraphics[width=0.6\textwidth]{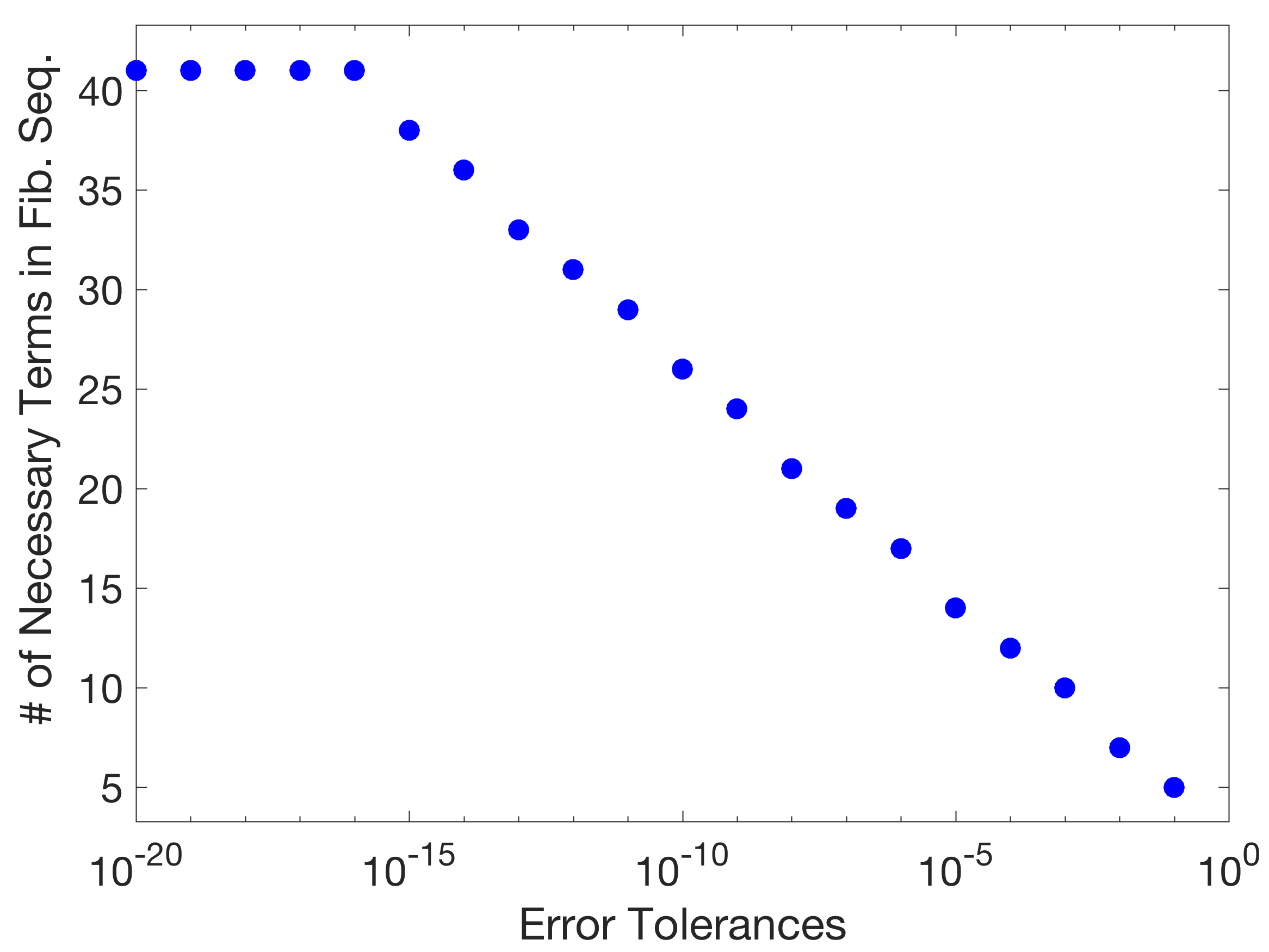}
    \caption{The number of terms \textit{necessary} to achieve a specific error tolerance vs. the error tolerance imposed.}
    \label{fig:GoldenRatio_numTerms}
\end{figure}

Another way we could have approached this problem is asking how many terms do we need to achieve a specific accuracy. This is subtly different than the previous question, where we calculated subsequent terms in the sequence and checked the error in each successive approximation. Now, for a particular \textit{error tolerance} specified, we will keep including additional terms in the Fibonacci Sequence until our approximation is within the error tolerance. %\textcolor{blue}{A script to do just this is found in\ldots}. 
Figure \ref{fig:GoldenRatio_numTerms} illustrates this relationship.   

Figure \ref{fig:GoldenRatio_numTerms} suggests that for error tolerances smaller than $10^{-15}$, we still only need $41$ terms, however, this is again a computational artifact due to constraints of double precision (\textit{64-bit}) accuracy. 

In general, a plot such as in Figure \ref{fig:GoldenRatio_numTerms} would be ideal for scientific computing applications. That is, we specify an error tolerance and we can figure out how much \textit{resolution} (in this case, number of terms in the Fibonacci Sequence) we need to achieve it. Unfortunately, in practice this tends to not be the case. What allowed us to do that here was the fact that we are able to compute successive approximations to $\phi_n$ extremely fast! In obtaining the data in Figure \ref{fig:GoldenRatio_numTerms} we are actually performing a lot of extra \textit{floating point operations} (operations) to answer the question of how many terms we need to achieve a certain accuracy. 

In research applications, many codes will not run anywhere near this quickly, so rather than asking ``how much resolution (\# of terms) do I need for a certain accuracy?" we are forced to try certain resolutions and then check the accuracy, as in Figure \ref{fig:GoldenRatio_Conv}. Moreover, in most cases since the true solution is not known, the best one can do is approximate the error against more highly resolved cases. This is precisely what we will illustrate in Sections \ref{sec:trap_rule} and \ref{jellyfish} when we approximate the value of a non-trivial definite integral and introduce a jellyfish locomotion model. Furthermore in Section \ref{jellyfish} we will see that the simulations require non-trivial computational time to run (from minutes to weeks!) and so we must be mindful of not performing excess computations than are essential. 

To summarize this section we have seen that (1) convergence plots compare the error against the resolution used (here the number of Fibonacci terms), (2) additional terms in the Fibonacci Sequence make approximations more accurate, (3) the Golden Ratio Sequence exhibits geometric convergence, and (4) if accuracies of $10^{-15}$ are obtained, they are the limit of standard \textit{64-bit} double-precision computers. 

Scripts to perform the computations in this section, as well as an example class activity are found in {\em Supplemental/Golden$\_$Ratio/}. In the next section, we investigate how the convergence properties of a numerical algorithm may change depending on the problem that it is applied to. 

%-we see that finding the value of the golden ratio on a computer in this manner seems almost instantaneous, we will show jellyfish example in which speed of simulation is a determining factor in what error to live with. 

%
% LEARNING GOALS
%
%1. See convergence plot in action
%2. See additional terms make more accurate
%3. Computers are only accurate to double precision accuracy

%%%%%%%%%%%%%%%%%%%%%%%%%%%%%%%%%%%%%%%%%%%%%%%%%%%%%%%%%%%%%%%%%%%%%%%
%
% TRAP RULE
%
%%%%%%%%%%%%%%%%%%%%%%%%%%%%%%%%%%%%%%%%%%%%%%%%%%%%%%%%%%%%%%%%%%%%%%%

%%%%%%%%%%%%%%%%%%%%%%%%%%%%%%%%%%%%%%%%%%%%%%%%
%
%
% TRAP RULE
%
%
%%%%%%%%%%%%%%%%%%%%%%%%%%%%%%%%%%%%%%%%%%%%%%%%

\section{Composite Trapezoid Rule Convergence}
\label{sec:trap_rule}

In this Section we illustrate that the convergence rate of a numerical scheme may depend on the problem to which it is applied  \cite{Eggert:1989,Trefethen:2014a,Trefethen:2014b}.

\begin{figure}[H]
    %\centering
    \centering
    \includegraphics[width=0.65\textwidth]{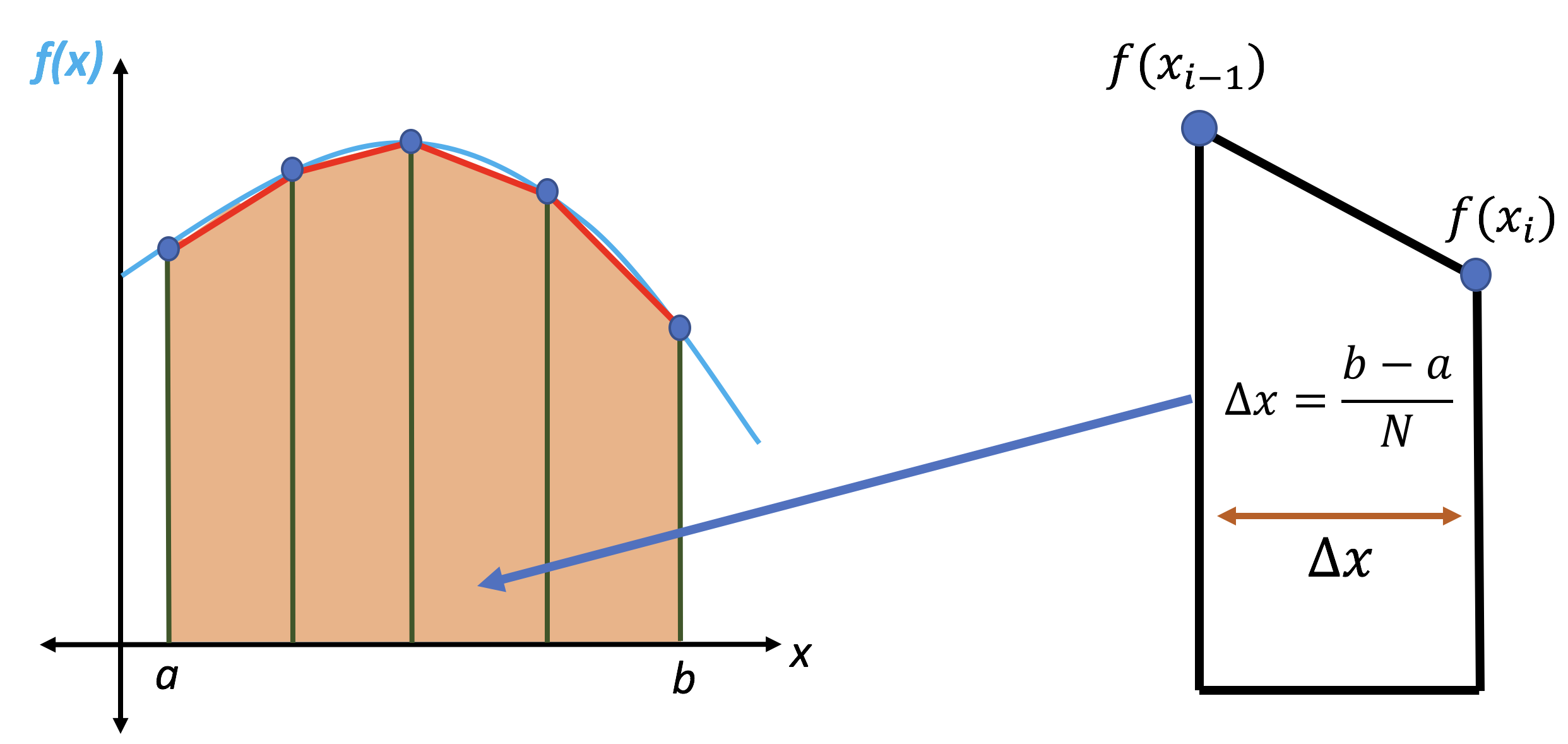}
    \caption{Image illustrating the the composite trapezoid rule.}
    \label{fig:TrapRule}
\end{figure}

Typically during second semester Calculus students are introduced to the trapezoid rule, in particular the \textit{composite} trapezoid rule for approximating definite integrals. The motivation usually stems from the integrand not having have a closed form anti-derivative. For an integral such as
\begin{equation}
    \label{eq:int} I = \int_{a}^{b} f(x) \ dx
\end{equation}
recall that its approximation using the composite trapezoid rule using $N$ partitions (or subintervals) is given by 

$$I_{N} = \frac{b-a}{2N}\Big[ f(x_0=a) + 2f( x_1 ) + 2f(x_2) + \ldots 2f(x_{N-1}) + f(x_N=b) \Big],$$

where $x_j = a + j\Delta x$ for $j=0,1,2\ldots,N$ and $\Delta x$ is the width of the base of each trapezoid, see Figure \ref{fig:TrapRule}.

Beyond simply approximating an integral using $5$, $10$, or (hopefully not more than) $20$ partitions by hand using the trapezoid rule, students in Calculus are taught a formula to find an error bound for this approximation using $N$ partitions.
%If $E_T$ is the error associated with using the composite trapezoid rule, 
If $|E_T|=|I-I_N|$ is the {\em absolute error} associated with using the composite trapezoid rule, then the error is bounded by
\begin{equation}
    \label{eq:TrapErrBound} |E_T|\leq \frac{K(b-a)^3}{12N^2}.
\end{equation}
In this error bound, $a$ and $b$ are the integration bounds and $N$ is the number of partitions; the only unknown quantity is $K$. Since we are told that the approximation error is \textit{no more} than the right hand side of (\ref{eq:TrapErrBound}), one can motivate the quantity $K$ to have a clear flavor of optimization to it in one way or another. Since we never want to underestimate the error in an approximation, we seek the \textit{worst-case} error bound.%meaning that we could choose $K$ to be as large as possible for the problem. 

For the composite trapezoid rule we find that 
\begin{equation} 
K=\underset{a\leq x\leq b}{\max}\Big|f''(x) \Big|.\label{eq:k_value}
\end{equation} %\ \mbox{for all} \ a\leq x\leq b.$$

%That is to say, $K$ is the maximum of the second derivative of the integrand $f(x)$ over $a\leq x \leq b.$ 
%This choice of $K$ ensures that we will not underestimate the error. 
Students who have taken a numerical analysis class may recognize this error bound also, as they may have proved it using Taylor Series expansions \cite{Kincaid:2002,Burden:2014}.

%%%%%%%%%%%%%%%%%%%%%%%%%%%%%%%%%%%%
% TRAP RULE: NON-PERIODIC
%%%%%%%%%%%%%%%%%%%%%%%%%%%%%%%%%%%%

\begin{exmp}
\label{ex:trap_NP} Let's apply the composite trapezoid rule with uniform partition spacing to the following integral, $$I = \int_{0}^{1} \frac{(x^2+3)cos^2(2\pi x)}{ ( 1 + e^{sin(2\pi x) })^2}\ dx.$$

This does not look like an integral where we could analytically compute its anti-derivative (nor do we wish to try here). Applying the composite trapezoid for a variety of partitions, $N$, we will see that more partitions leads to higher accuracy, as (\ref{eq:TrapErrBound}) suggests. How accurate is it? This is interesting - we do not know the exact value of the integral; however, from \eqref{eq:k_value} it follows that we can approximate its solution with \textit{an appropriate} (``very large") number of partitions and use that value as though it were the exact solution (as we used $\hat{\phi}$ as the ``exact'' value of the Golden ratio in Section \ref{sec:golden_ratio}). For $N=10,000,000$, the approximation is $$I_{10000000}=0.455122322888408.$$

\begin{figure}[H]
    %\centering
    \centering
    \includegraphics[width=0.65\textwidth]{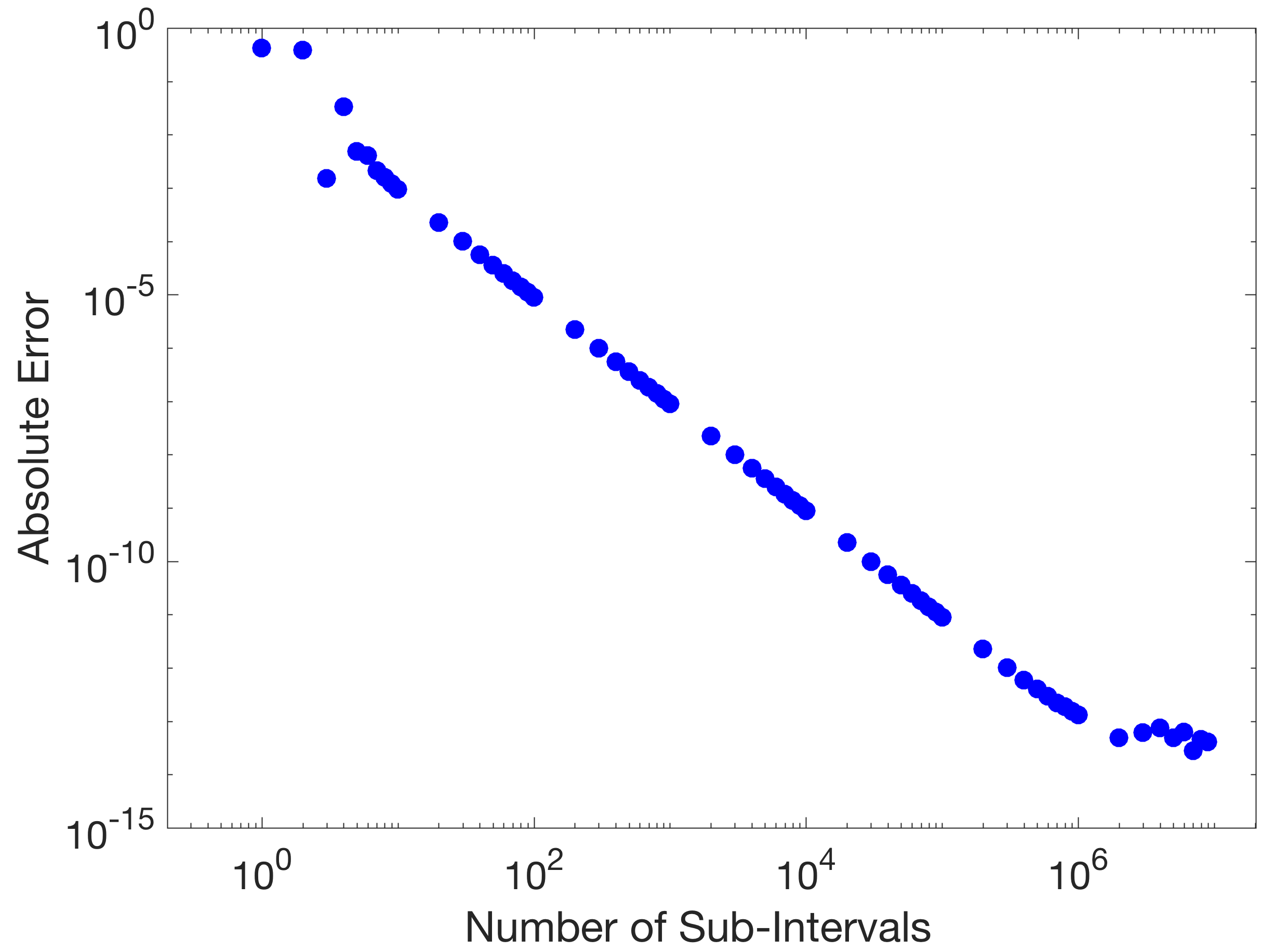}
    \caption{Plot depicting how the absolute error decreases as the number of partitions, e.g., number of trapezoids, increases in the integration domain.}
    \label{fig:TrapRuleErrNP}
\end{figure}

Figure \ref{fig:TrapRuleErrNP} showcases how the error decreases as more partitions are added, when compared to the ``true" solution that we approximated above using $10,000,000$ partitions. It takes over $1$ million partitions until we reach machine precision in this case! Furthermore, this figure illustrates a linear relationship between the logarithm of absolute error and logarithm of number of partitions used. Hence as in Section \ref{intro} we see the following error relationship

$$\mbox{Absolute Error}\sim N^{\mbox{slope}},$$

where the slope is negative. Computing the slope of the line once the resolution is high enough, e.g., there are enough partitions $\sim 10^2$, we find that $\mbox{slope}\approx -1.978831 \approx -2$, in agreement with (\ref{eq:TrapErrBound}). 

\end{exmp}
Next we will perform a similar calculation, but on a definite integral whose integrand is \underline{periodic} on the integration domain.

%\textcolor{red}{General question: can one make sense of whether the machine precision threshold occurs due to (i) finite storage of numbers in the computer versus (ii) after many calculations there is accumulation of errors which add up?}

%%%%%%%%%%%%%%%%%%%%%%%%%%%%%%%%%%%%
% TRAP RULE: PERIODIC
%%%%%%%%%%%%%%%%%%%%%%%%%%%%%%%%%%%%
\begin{exmp}
\label{ex:trap_P} We will now apply the composite trapezoid rule with uniform partition spacing to the following definite integral, $$I = \int_{0}^{1} \frac{cos^2(2\pi x)}{ ( 1 + e^{sin(2\pi x) })^2}\ dx.$$

Note that this integral is very similar to the previous integral, except with one special property - the integrand is now \underline{periodic} on the integration domain. Since we still do not know the \textit{exact} value of this integral, we again use 10,000,000 partitions in a composite trapezoid rule, to determine $$I_{10000000} = 0.132214293037990.$$ 
%we will use the same trick to obtain a very accurate approximation - use \textit{enough} partitions. Using 10,000,000 partitions again in a composite trapezoid rule, we find our approximation to be $$I_{10000000} = 0.132214293037990.$$ 

Figure \ref{fig:TrapRuleErrP}a gives a log-log convergence plot for the absolute error vs. number of sub-interval partitions; however, it is clear that there is a distinct difference between this plot and the convergence plot in Figure \ref{fig:TrapRuleErrNP}; the accuracy appears to achieve machine precision extremely quickly! In fact, it appears to achieve this accuracy for $N\sim10$. Figure \ref{fig:TrapRuleErrP}b highlights this accelerated convergence, by showing the same data, but on a semi-logarithmic plot with fewer partitions. 

\begin{figure}[H]
    %\centering
    \centering
    \includegraphics[width=0.99\textwidth]{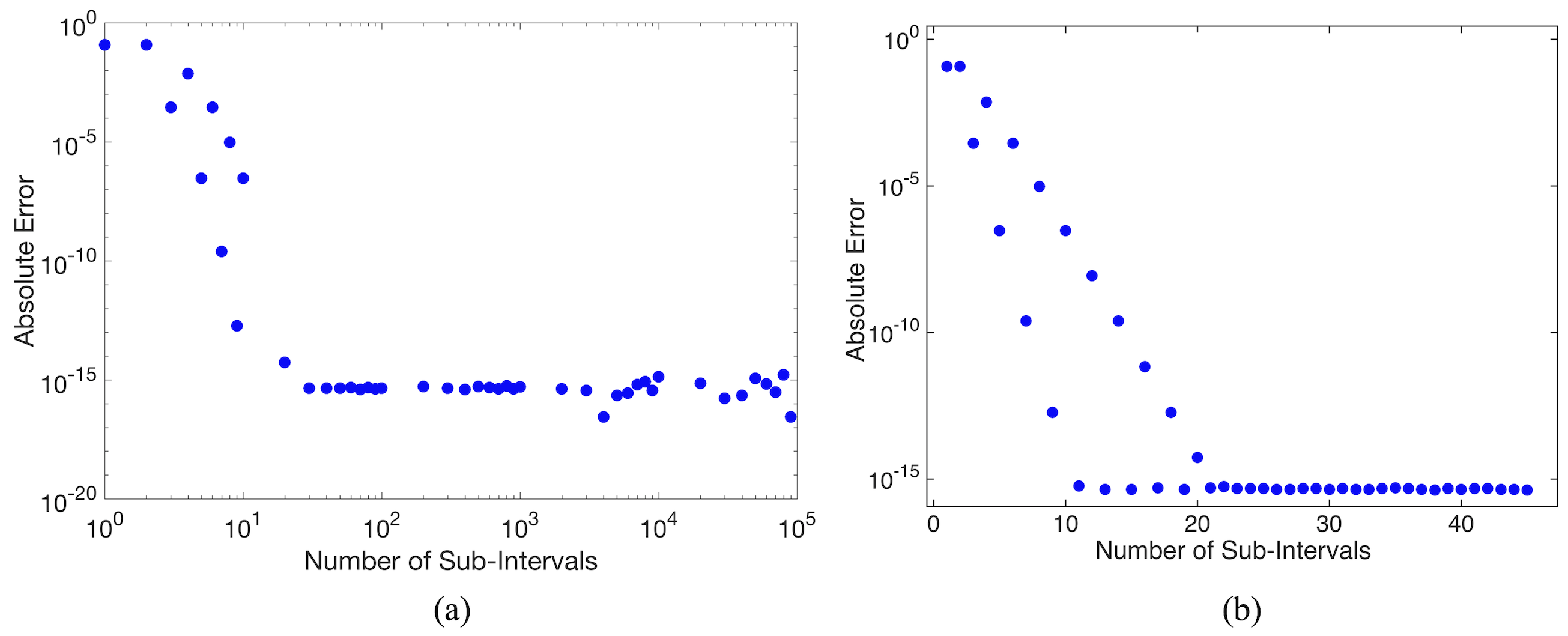}
    \caption{Illustrating the convergence rate of the composite trapezoid rule for an example of a definite integral with an integrand that is periodic on the integration domain.}
    \label{fig:TrapRuleErrP}
\end{figure}

Another aspect that is particularly interesting is that it appears that before achieving machine precision accuracy, there are two different convergence rates, corresponding to whether $N$ is odd or even. Figure \ref{fig:TrapRuleErrP}b says that for the approximations using an odd number of partitions, they converge quicker towards the ``exact" solution. However, whether odd or even, in this case the composite trapezoid rule achieved machine precision only using $21$ (or less) sub-intervals! Recall that in the Example \ref{ex:trap_NP} it took almost 1 million partitions. That is a significant difference! Moreover, from witnessing a linear relationship between the logarithm of the absolute error and number of sub-intervals in Figure \ref{fig:TrapRuleErrP})b, we see that the error decreases \textit{geometrically} as in Section \ref{sec:golden_ratio}. Appendix \ref{app:TrapRuleExpConv} provides justification for this accelerated convergence rate. For a deeper study on when the trapezoid rule converges geometrically, see \cite{Trefethen:2014a,Trefethen:2014b}.

\end{exmp} % END PERIODIC CASE CONVERGENCE

Comparing the convergence rates from Examples \ref{ex:trap_NP} and \ref{ex:trap_P}, we have just seen is that the composite trapezoid rule can exhibit significantly different convergence rates depending on what definite integral it is applied to.  % Although this is not a new phenomenon \cite{Eggert:1989,Trefethen:2014a,Trefethen:2014b} it illustrates that the convergence rate of a numerical scheme is not static among all problems it could be used for. 

In these examples if we only cared about obtaining a particular accuracy, it may have been enough to simply use the number of partitions as $N=10,000,000$ and not looked back. It is clear that one case almost warranted this much \textit{resolution}, while in another case it was extreme overkill. You might be thinking \textit{well, either way it computed the definite integral accurately for $N=10,000,000$}. However, when running the script we see that as $N$ increases, it takes longer to compute the integral approximation. If computing a definite integral was simply one part, e.g., module, of many in a larger numerical scheme, the difference between finding an integral approximation that takes $0.001s$ and $1.0s$ could be a significant difference if, say, you need to integrate thousands of definite integrals when running the algorithm. To that end, as we've seen from the two different convergence plots, this may be akin to asking ourselves how much error we (our application) can tolerate, e.g., still obtaining physical results with only $10^{-6}$ accuracy vs. $10^{-12}$ accuracy, as in Figure \ref{fig:TrapRuleErrNP}.

The script that will run these examples as well as an example of a class activity is found in \textit{Supplemental/Trapezoid$\_$Rule/}. To summarize this section, we have seen (1) the composite Trapezoid Rule can exhibit geometric convergence if applied to a definite integral with periodic integrand, (2) unexpected and complex convergence properties arise from even seemingly basic numerical schemes,  and (3) that the convergence rate of a numerical scheme can significantly change depending on the problem to which it is applied.
 
Next in Section \ref{jellyfish} we will investigate a contemporary computational model of jellyfish locomotion where all these ideas will be extensively explored: no longer will we be able to simply increase the grid resolution indefinitely, for doing so will necessitate increases of computational time from mere minutes to days to weeks or even months!

%%%%%%%%%%%%%%%%%%%%%%%%%%%%%%%%%%%%%%%%%%%%%%%%%%%%%%%%%%%%%%%%%%%%%%%
%
% JELLYFISH 
%
%%%%%%%%%%%%%%%%%%%%%%%%%%%%%%%%%%%%%%%%%%%%%%%%%%%%%%%%%%%%%%%%%%%%%%%

%%%%%%%%%%%%%%%%%%%%%%%%%%%%%%%%%%%%%%%%%%%%%%%%
%
%
% JELLYFISH CONVERGENCE
%
%
%%%%%%%%%%%%%%%%%%%%%%%%%%%%%%%%%%%%%%%%%%%%%%%%

\section{Jellyfish Locomotion: Convergence and Speed!}
\label{jellyfish}

As previously seen in Sections \ref{sec:golden_ratio} and \ref{sec:trap_rule} a convergence plot is useful for figuring out how many terms (how much \textit{resolution}) is necessary to achieve a certain accuracy. In this section we will dive head first down that rabbit hole for an example involving jellyfish locomotion; however, there will be three distinctions: (1) We do not have an analytic result to compare our numerical results nor one that be trivially obtained, (2) we will not even have a straight-forward metric in which to compute the error, and (3) increasing the resolution, while decreasing the error, results in significant increases in the computational time required to run a simulation, mandating mindfulness when changing the resolution (here given by spatial grid-steps).

\begin{figure}[H]
    %\centering
    \centering
    \includegraphics[width=0.99\textwidth]{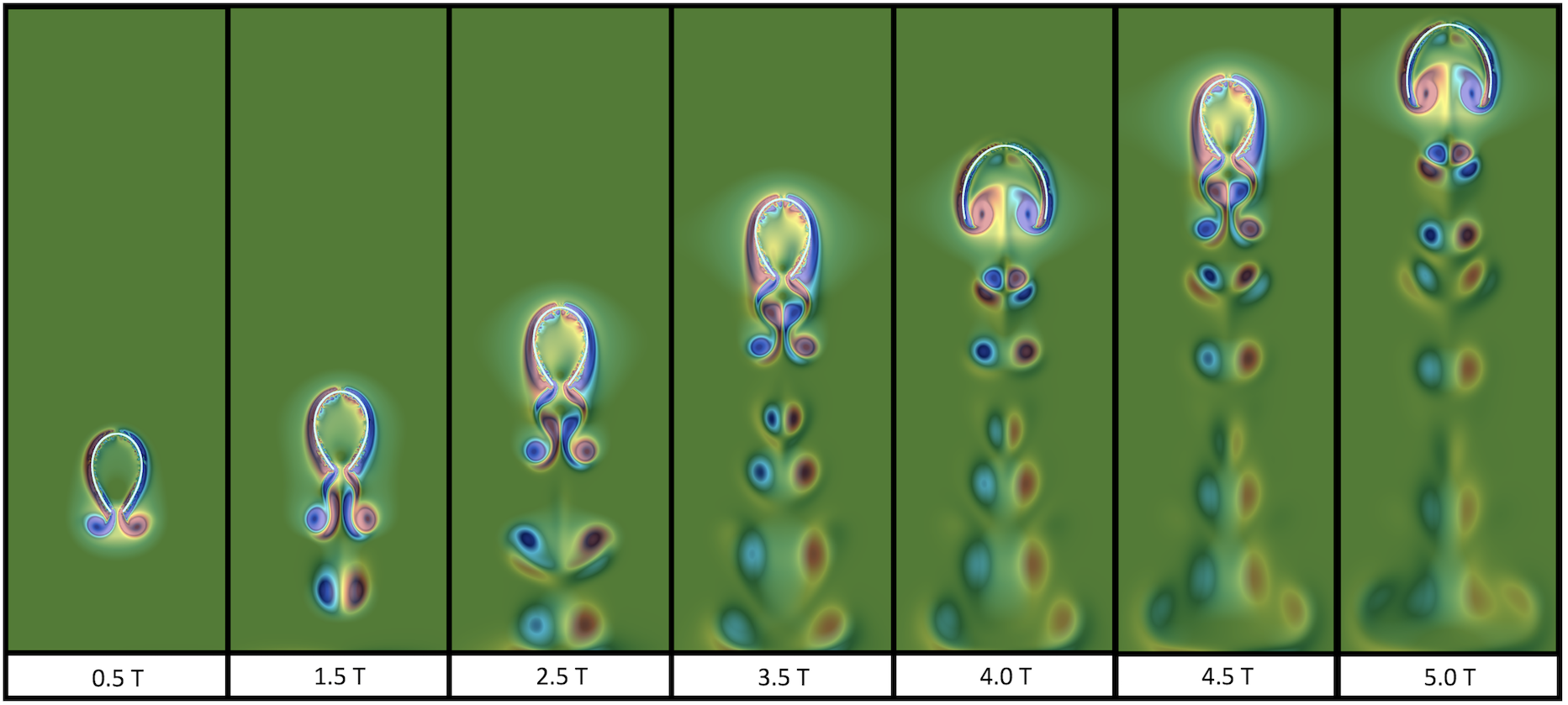}
    \caption{A simulation of a single jellyfish swimming at $Re=150$. The background colormap represents vorticity and $T$ refers to the number of periods of jellyfish bell contractions.}
    \label{fig:Single_Jelly}
\end{figure}

This model of jellyfish locomotion can be found in the \textit{IB2d} example folder \textit{Example$\_$Jellyfish$\_$Swimming/Hoover$\_$Jellyfish/} as well in the Supplemental materials, \textit{Supplemental/Jellyfish/Simulation$\_$Skeletons/}. We note that the code is initialized on a rectangular domain with resolution $512\times 512$, at Reynolds Number $Re=150$, as in the model to which it is based \cite{Hoover:2015}. For more information regarding Reynolds Number see Appendix \ref{app:Re}. Snapshots of this jellyfish swimming are found in Figure \ref{fig:Single_Jelly}; it is clear that the jellyfish is able to produce enough forward thrust to propel itself forward, that is, swim forward at this resolution and $Re$. 

An unfortunate, yet inescapable, aspect about the simulation shown in Figure \ref{fig:Single_Jelly} is that it takes roughly 8 hours to run on a personal computer (iMac with 6 GB 2400 MHz DDR4 memory and 3.6 GHz Intel Core i7 processor). Why does it take so long? Fluid-structure interaction simulations are notorious for being computationally expensive and extensive effort has placed into reducing simulation time \cite{McQueen1997,Roma:1999,GriffithThesis:2005, Griffith:2007,BGriffithIBAMR}. Rather than jump into the intricacies of such techniques, such as parallelization \cite{McQueen1997} or adaptive mesh refinement strategies \cite{Roma:1999,Griffith:2007}, we will focus our efforts on simply comparing computational time for increasing grid resolution and accuracy. %Recall that with increasing resolutions comes increased accuracy, but as we will see it comes at the price of significantly increasing the computational time required for a simulation to run. 

We are now entering a practical domain in scientific computing; how much accuracy is required and how can I achieve that in a \textit{timely} manner? That is, if it takes one week for a simulation to give us 8 decimals of accuracy, but it only takes one day to give us 6 decimal accuracy, is the extra six days of computational time worth those extra two decimals of accuracy? This is a philosophical question and of course, depending on the application, the answer could very well be ``surely no" or ``absolutely". %Here we give a sense of what computational research is like, including the very practical trade-offs between accuracy and computational expense. 
These ideas will be, you guessed it, quantified through the use of convergence plots. 

Before we can discuss convergence plots, we need to explicitly define what we mean by error in a fluid-structure interaction simulation. It has been defined through the differences between either fluid (Eulerian) or immersed structure (Lagrangian) data, e.g., position of the immersed body \cite{Griffith:2005,Mori:2008,Griffith:2007}, fluid velocity \cite{Lai:2000,Griffith:2005,Mori:2008,Griffith:2007}, pressure \cite{Griffith:2005,Griffith:2007}, or forces on the immersed body \cite{Lai:2000,BattistaIB2d:2016}. These are the quantities used in determining error, however if we are computing the error for those quantities at different grid resolutions, what are we using as the \textit{true} or \textit{exact} answer in which to compare them? This is a very important question in research grade numerical analysis. The answer is actually a practical one - we compare those quantities to a simulation with very high accuracy, usually the highest accuracy that is realistically possible in a \textit{timely manner}. %Again with the timely manner phrasing. 
In this case, an \textit{untimely manner} describes a simulation that takes significantly longer than the simulations that you actually want to use for data collection. 

For our purposes here, we will compare jellyfish locomotion on $N\times N$ grids where $N\in\{32,48,64,96,128,256,384,512,768,1024,1536,2048\}$. Recall that the simulation with $N=512$ took approximately 8 hours - just imagine how long a simulation with $N=1024$ or $N=2048$ takes? We can give a crude estimate, simply by noting that if we go from $512\rightarrow1024$ that is a factor of $2$ and since we are in two-dimensions, that gives a computational time scaling of $2^2$. Therefore going from $512\rightarrow1024$ will likely take \textit{4x} as long as the $512$ simulation! Then going from $512\rightarrow2048$ is a factor of $4$ which gives a time scaling of $4^2=16$, estimating it would take approximately \textit{16x} as long! In general this computational expense scaling can be written as 

\begin{equation}
    \label{eq:COD} \mbox{estimated comp. time scaling} = (\mbox{resolution factor increase})^{\mbox{dimension}}.
\end{equation}

\begin{figure}[H]
    %\centering
    \centering
    \includegraphics[width=0.65\textwidth]{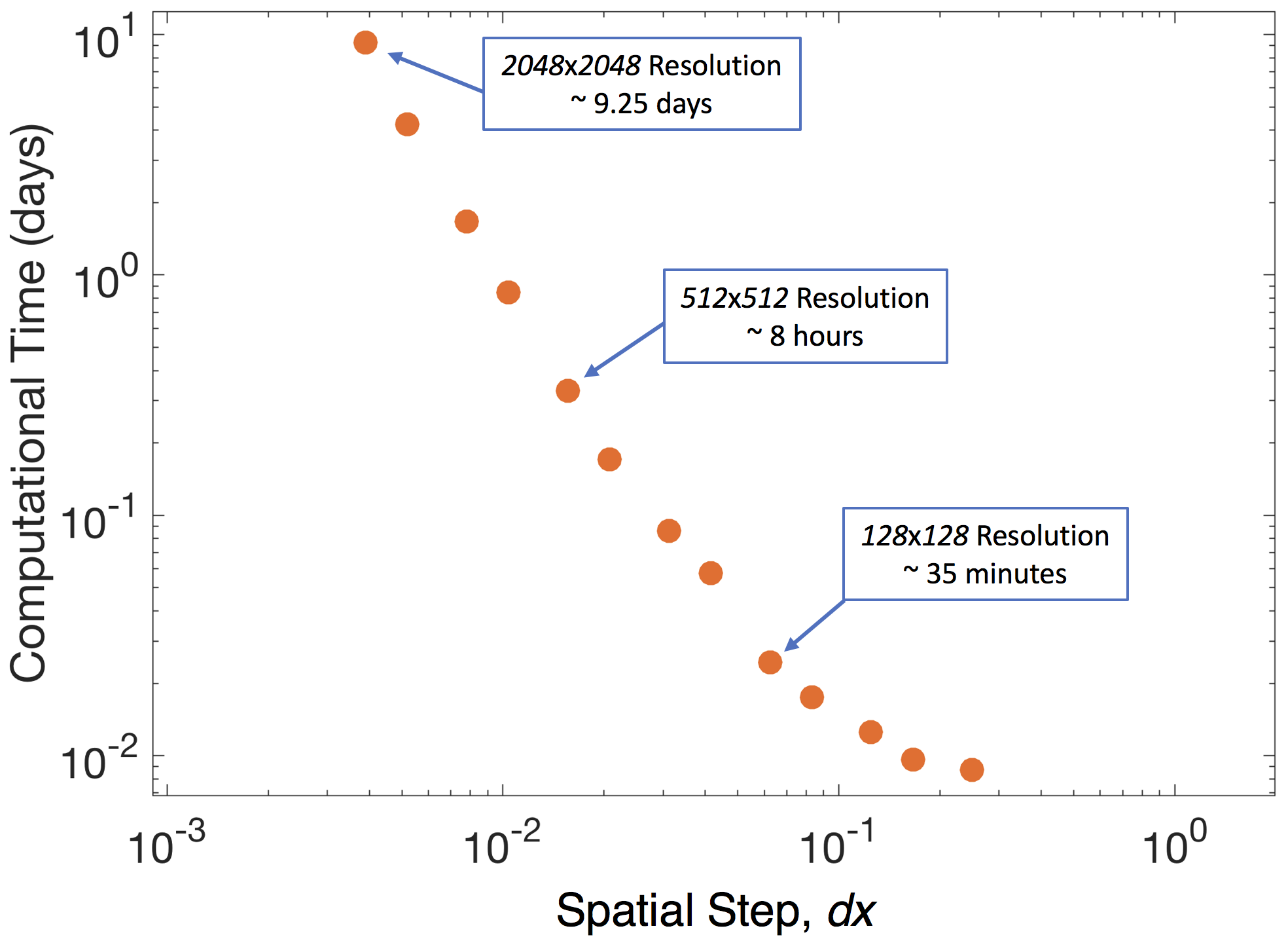}
    \caption{Computational time required to for $5$ jellyfish bell contractions. As the computational grid gets more refined ($dx$ gets smaller) so accuracy increases, computational time also increases geometrically.}
    \label{fig:Jelly_CompTime}
\end{figure}

If we were performing simulations in $3D$, the scaling becomes even more extreme! This is known as the \textit{curse of dimensionality} \cite{Chen:2009}. Luckily we are only working in two dimensions here. For simulations consisting of $5$ jellyfish pulses, with time-step, $dt =1e-5$, a plot of the computational time required versus grid resolution is given in Figure \ref{fig:Jelly_CompTime}. These simulations were run on an iMac with a 3.6 GHz quad-core $7^{th}$ generation Intel processor with 16 GB 2400 Mhz DDR4 memory. There is a large difference between the $512$ and $2048$ resolution cases. 

%%%%%%%%%%%%%%%%%%%%%%%%%%%%%%%%%%%%%%%%%%%%%%%%
%
% QUALITATIVE JELLYFISH CONVERGENCE
%
%%%%%%%%%%%%%%%%%%%%%%%%%%%%%%%%%%%%%%%%%%%%%%%%

\subsection{Qualitative Convergence}
\label{sec:qualitative_jellyfish}

\begin{figure}[H]
    %\centering
    \centering
    \includegraphics[width=0.99\textwidth]{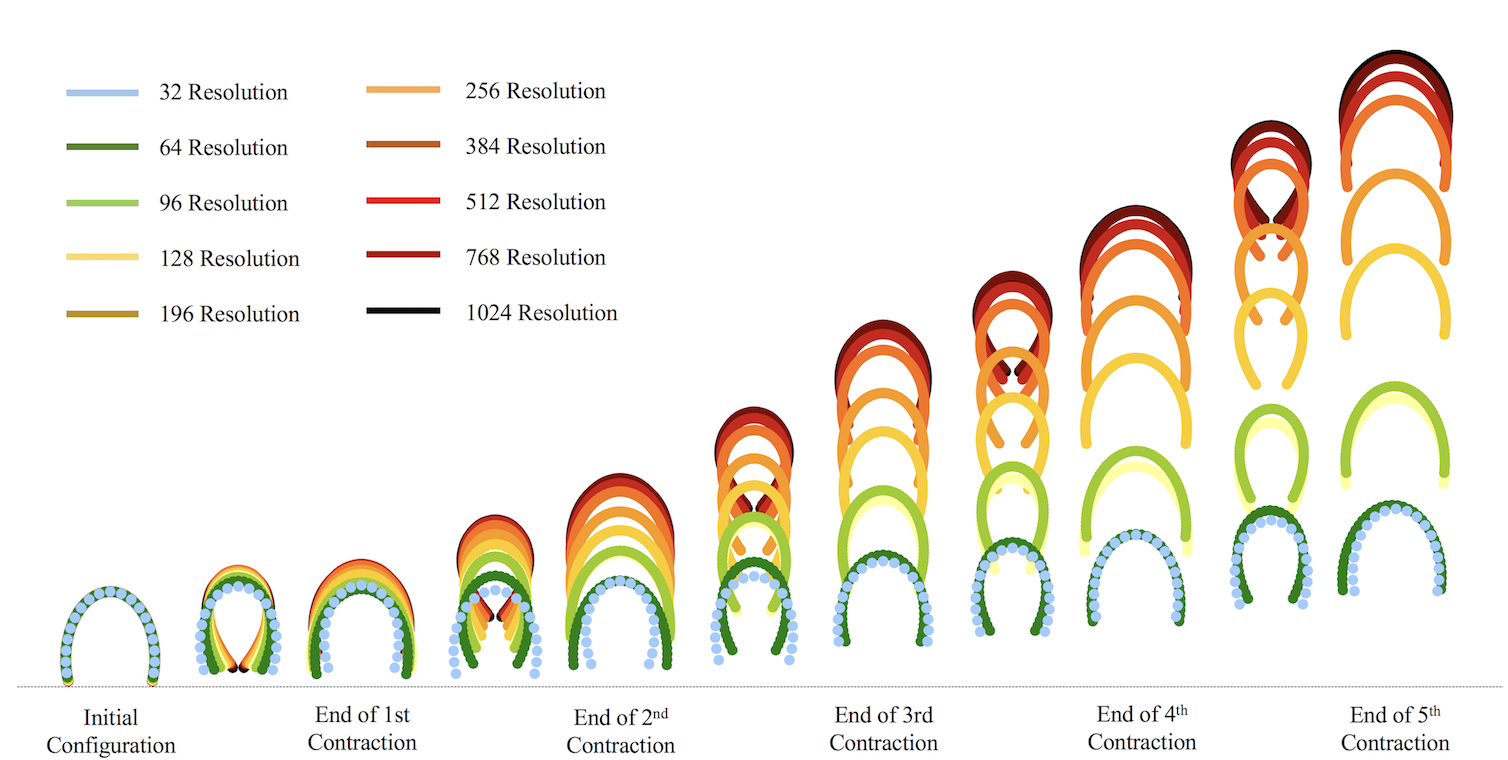}
    \caption{Jellyfish swimming performance at varying grid resolutions. Better forward swimming performance is achieved at higher grid resolutions.}
    \label{fig:Jelly_Conv_Sims}
\end{figure}

As suggested previously, there are many ways we can quantify error in these simulations. Before explicitly computing the error, let's simply observe if we can qualitatively see differences between simulations with different accuracies by visualization. Figure \ref{fig:Jelly_Conv_Sims} shows multiple simulations overlain at various uniform time points during the simulation. In each of these cases, the jellyfish uses the same physical mechanisms for propulsion, that is, the contraction forces and material properties of the jellyfish itself are scaled appropriately. Note that these simulations were carried out at $Re=150$ (see Appendix \ref{app:Re}).

It is clear that better forward swimming performance is achieved for higher grid resolutions, but by $512\times 512$ resolution the simulations become quite similar. It is common that in published manuscripts related to jellyfish swimming that grid resolutions are in the realm of $dx\sim0.015625$ on domains of length 8, e.g., $512\times512$ spatial resolutions (\cite{Hershlag:2011,Alben:2013,Hoover:2015,Hoover:2017,Hoover:2019}). 

%%%%%%%%%%%%%%%%%%%%%%%%%%%%%%%%%%%%%%%%%%%%%%%%
%
% QUANTITATIVE JELLYFISH CONVERGENCE
%
%%%%%%%%%%%%%%%%%%%%%%%%%%%%%%%%%%%%%%%%%%%%%%%%

\subsection{Quantitative Convergence}
\label{sec:quanitative_jellyfish}

In this subsection we will define various metrics for the error in these jellyfish simulations. There will be errors defined on the Lagrangian domain (e.g., the jellyfish) and errors defined on the Eulerian grid (e.g., the fluid grid). 

%%%%%%%%%%%%%%%%%%%%%%%%%%%%%%%%%%%%%%%%%%%%%%%%
%
% QUANTITATIVE JELLYFISH CONVERGENCE: on jelly
%
%%%%%%%%%%%%%%%%%%%%%%%%%%%%%%%%%%%%%%%%%%%%%%%%

\subsubsection{Error on Lagrangian Structure (jellyfish)}
\label{sec:lagrangian_err}

\begin{figure}[H]
    %\centering
    \centering
    \includegraphics[width=0.65\textwidth]{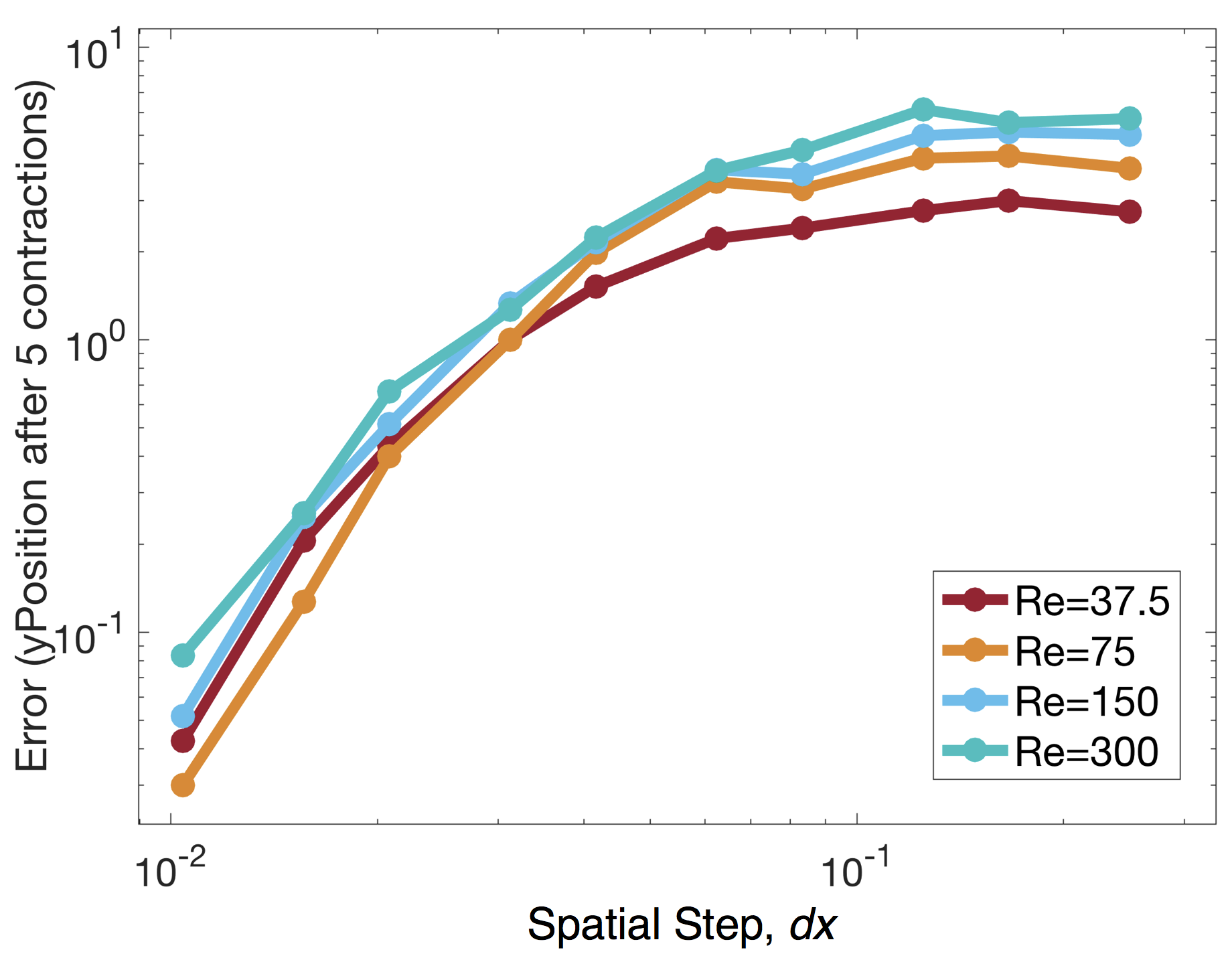}
    \caption{Illustrates the error defined by the difference between the top of the bell in the highest resolved case of $dx=0.0078125$ ($1024\times1024$ resolution) and less resolved grids for $Re=\{37.5,75,150,300\}.$}
    \label{fig:Jelly_Conv_Sims_yPos}
\end{figure}

To spin off of Section \ref{sec:qualitative_jellyfish}, we could begin quantifying convergence by defining the absolute error to be the spatial difference between the top of the bell with resolution $1024\times1024$ to less resolved cases, e.g., 
\begin{equation}
    \label{eq:jelly_yPos_err} err_{\mbox{y}} = | y_{1024}^5 - y_{j}^5|, 
\end{equation}
where $j=32,48,\ldots,768$, and $y_j^5$ is $y$ position of the center of the bell after $5$ contraction cycles for a particular resolution $j$. As seen in Figure \ref{fig:Jelly_Conv_Sims}, better forward swimming performance was associated with more highly resolved simulations, that is, these model jellyfish do not swim effectively at low grid resolutions. Figure \ref{fig:Jelly_Conv_Sims_yPos} gives such error vs. grid resolution for different batches of simulations with uniform $Re$, $Re=\{37.5,75,150,300\}.$ It illustrates that it is not until the spatial step $dx\sim0.0625$ corresponding to the $128\times128$ resolution, that the error begins to significantly decrease in all cases of $Re$ given. Moreover, it is apparent that the amount of error in bell position is also a function of the Reynolds Number, $Re$, being considered.

%Beyond simply defining the error to be the difference between the positions of the jellyfish bells,

\begin{figure}[H]
    %\centering
    \centering
    \includegraphics[width=0.65\textwidth]{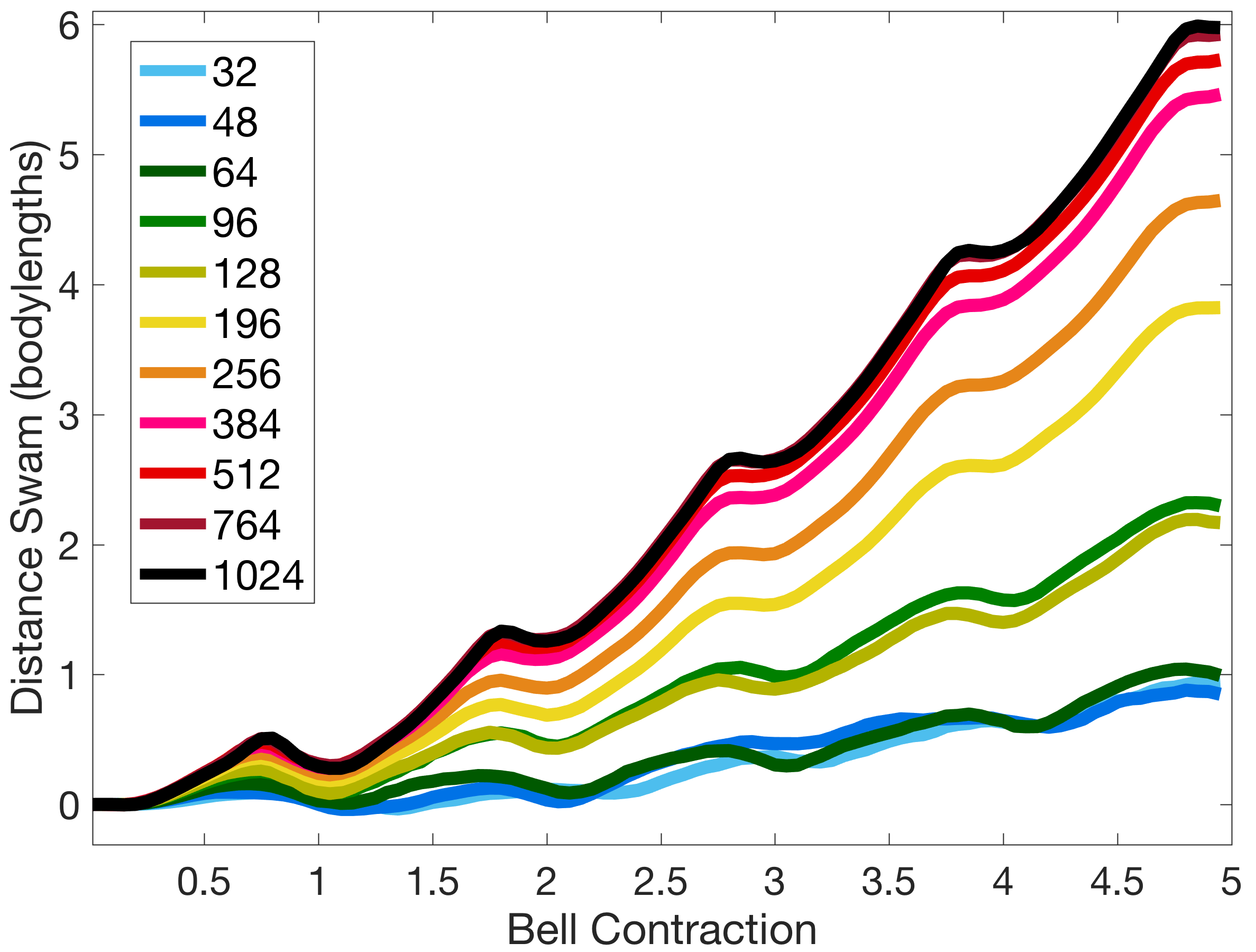}
    \caption{Gives the distance swam in jellyfish in bodylengths vs. bell contraction cycle. The higher resolution cases demonstrate better forward swimming performance.}
    \label{fig:Jelly_Conv_Sims_DistanceVsTime}
\end{figure}

We investigate the cause of this difference in positions by studying the error in the jellyfish velocities over time. To that end, we quantify an error associated with the swimming speed of the jellyfish after 5 bell contraction cycles. To compute the jellyfish swimming speed we first notice that the distance swam vs. bell contraction cycle appears close to linear, see Figure \ref{fig:Jelly_Conv_Sims_DistanceVsTime}. From Figure \ref{fig:Jelly_Conv_Sims_DistanceVsTime}, we find the best fit line through the data for each grid resolution tested during the last two bell contraction cycles and consider the slope of the line to be the average swimming speed for each case.

\begin{figure}[H]
    %\centering
    \centering
    \includegraphics[width=0.99\textwidth]{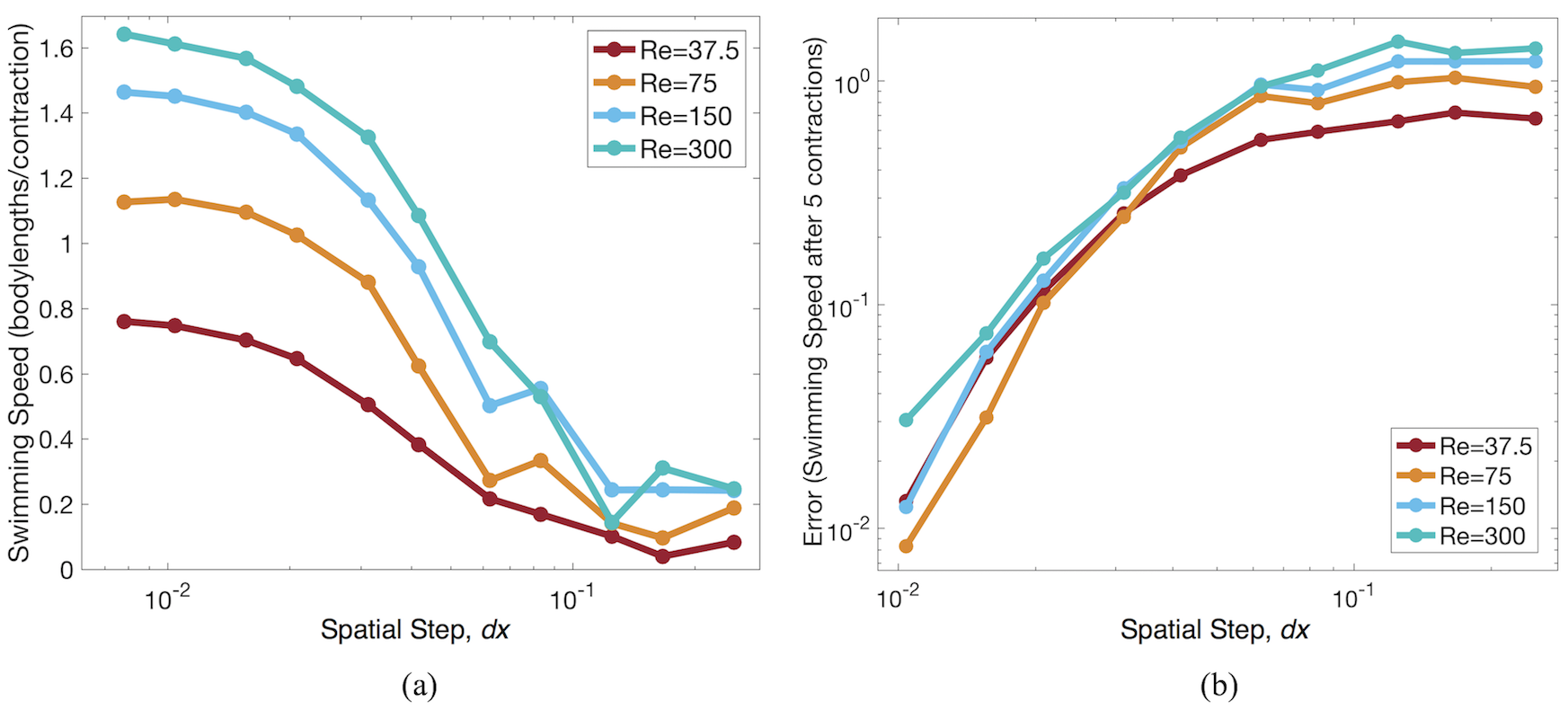}
    \caption{(a) The jellyfish's swimming speed after 5 contractions, averaged over the last two bell contraction cycles vs. grid resolution and (b) the absolute error given by the difference between the swimming speed at the highest resolution ($1024\times1024$) and less resolved grids.}
    \label{fig:Jelly_Conv_Sims_Speed}
\end{figure}

As the grid resolution increases, e.g., the spatial step, $dx$, decreases, the average swimming speeds begin to converge to the same value, see Figure \ref{fig:Jelly_Conv_Sims_Speed}a. The absolute swimming speed error is defined to be
\begin{equation}
    \label{eq:jelly_swim_err} err_{\mbox{speed}} = | m_{1024} - m_{j}|, 
\end{equation}
where $j=32,48,\ldots,768$, and $m_j$ is the slope of the best fit line computed from Figure  \ref{fig:Jelly_Conv_Sims_DistanceVsTime} for a particular resolution $j$. Similar to the position of the bell after 5 contractions, the swimming speed error is also a function of the $Re$ considered. As before, the absolute error drops off significantly for resolutions higher than $128\times128$. 

The errors observed in the distance swam (Figure \ref{fig:Jelly_Conv_Sims_yPos}) and swimming speed (Figure \ref{fig:Jelly_Conv_Sims_Speed}) arise from discrepancies in upward thrust (vertical) force generation when the jellyfish contracts at different grid resolutions. To investigate further, we analyzed the upward thrust force (y-directed force) on the jellyfish bell throughout $5$ bell contractions. The absolute error and relative error were defined as follows:

\begin{align}
    \label{eq:yForce_AbsErr} err_{ABS_{\mbox{thrust}}} &= \left| \bar{f}_{y_{1024}} ds_{1024} - \bar{f}_{y_{j}} ds_{j} \right|,  \\
    \label{eq:yForce_RelErr} err_{REL_{\mbox{thrust}}} &= \frac{ \left| \bar{f}_{y_{1024}} ds_{1024} - \bar{f}_{y_{j}} ds_j \right| }{ \left| \bar{f}_{y_{1024}} ds_{1024} \right| } = \frac{ err_{ABS_{\mbox{thrust}}} }{ \left| \bar{f}_{y_{1024}} ds_{1024} \right| }, 
\end{align}

where $\bar{f}_{y_j}$ denotes the spatially-averaged thrust (vertical) force across the jellyfish bell and $ds_j$ gives the Lagrangian distance between successive nodes on the bell. In both cases the $j$ denotes the specific spatial resolution, e.g., $j=32,48,\ldots,768.$ Multiplying by $ds_{j}$ converts from force densities to total force. 

\begin{figure}[H]
    %\centering
    \centering
    \includegraphics[width=0.925\textwidth]{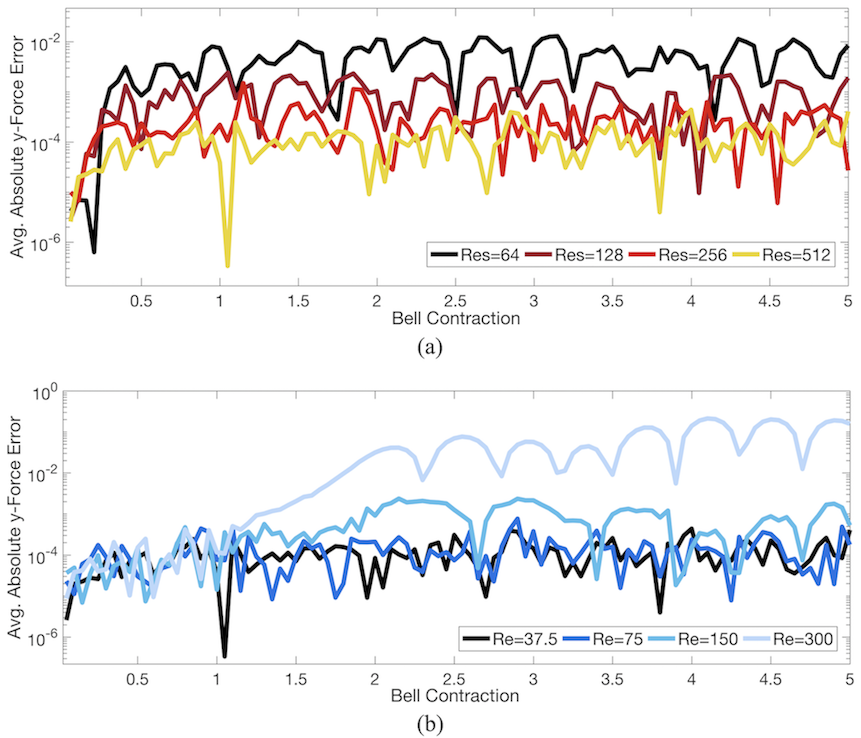}
    \caption{The absolute error given by the differences between the spatially-averaged thrust (vertical) force through $5$ bell contractions between the highest resolution ($1024\times1024$) and less resolved grids when (a) grid resolutions are varied and $Re=37.5$ and (b) the resolution is held at $512\times512$ and $Re$ is varied.}
    \label{fig:Jelly_Conv_Sims_Forces_Temporal}
\end{figure}

Figure \ref{fig:Jelly_Conv_Sims_Forces_Temporal} illustrates the resulting spatially-averaged vertical force through those contractions when (a) the grid resolution is varied while $Re=37.5$ and (b) the $Re$ is varied while grid resolution is held at $512\times512.$ From Figures \ref{fig:Jelly_Conv_Sims_Forces_Temporal}a and \ref{fig:Jelly_Conv_Sims_Forces_Temporal}b it is evident that as grid resolution increases, spatially-averaged vertical force error tends to decrease on average, as well that simulations using equivalent grid resolutions for different $Re$ do not lead to the same accuracies, respectively.

\begin{figure}[H]
    %\centering
    \centering
    \includegraphics[width=0.99\textwidth]{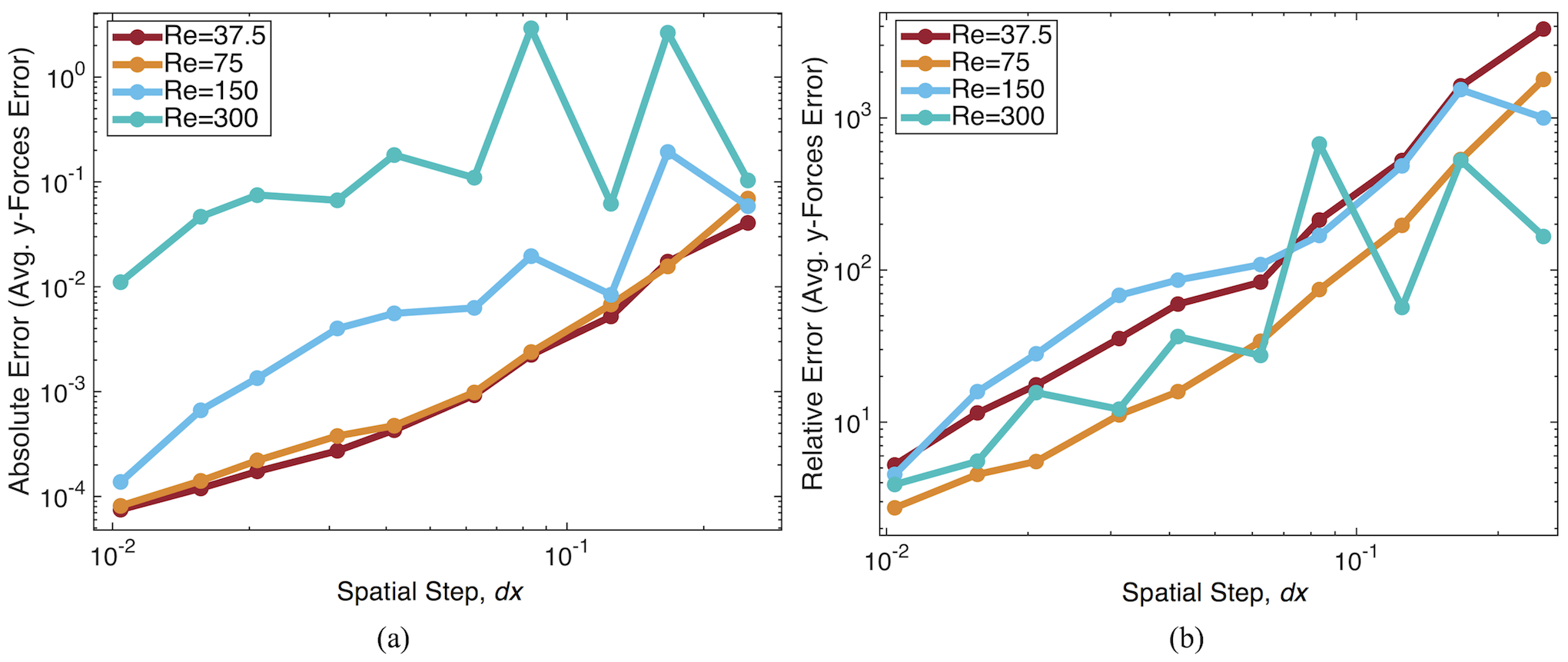}
    \caption{(a) The absolute error given by the difference between the spatially- and temporally-averaged thrust (vertical) force throughout the last swimming cycle between the highest resolution ($1024\times1024$) and less resolved grids. (b) The associated relative error computed analogously.}
    \label{fig:Jelly_Conv_Sims_Forces}
\end{figure}

Next we explored the spatially- and temporally-averaged upward thrust force (y-directed force) during the last bell contraction. The thrust force was averaged across the entire jellyfish bell for each grid resolution and Reynolds Number considered. Figure \ref{fig:Jelly_Conv_Sims_Forces} gives the absolute and relative errors, computed via Eqs.(\ref{eq:yForce_AbsErr})-(\ref{eq:yForce_RelErr}). Note that the relative errors are high because the assumed true value is small. The errors decay as grid resolution increases, e.g., the grid size decreases.

Furthermore, Figures \ref{fig:Jelly_Conv_Sims_yPos}, \ref{fig:Jelly_Conv_Sims_Speed}b, and \ref{fig:Jelly_Conv_Sims_Forces} illustrate an additional subtle point regarding mathematical models and their parameters - numerical errors may also depend on the parameter space explored, not only the time-step or spatial grid-size! Here for this jellyfish example using computational fluid dynamics - the numerical errors also depend on the Reynolds Number of the system! Although a solver was successfully used in one scenario, like $Re=37.5$, to achieve a certain accuracy, does not guarantee the same accuracy when applied to an almost identical system, but at higher $Re$, as in Figure \ref{fig:Jelly_Conv_Sims_Forces}a.

To summarize, in a holistic manner, the errors in upwards thrust force (Figure \ref{fig:Jelly_Conv_Sims_Forces}) for less resolved grids cause the jellyfish not to swim as fast (Figure \ref{fig:Jelly_Conv_Sims_Speed}) nor swim as far (Figure \ref{fig:Jelly_Conv_Sims_yPos}) as more resolved simulations. In short, less resolved grids lead to different locomotive patterns stemming from errors. However, as seen in Figure \ref{fig:Jelly_CompTime}, always pushing for higher resolutions may not be practical because of total computational time. 

To that note, one must inquire how much accuracy is required for a problem, e.g., validating that the jellyfish is capturing biologically relevant kinematics and/or swimming speeds at certain resolutions. Unfortunately, if obtaining the highest accuracy possible is the number one goal, computational models would be dramatically limited, as no model is ever complete or can describe nature exactly. In the words of G. E. Box, ``All models are wrong but some are useful." \cite{Box:1979}. For instance, the jellyfish model presented here only includes a $2D$ representation of its bell, one could include electrophysiology or porous tentacles or other complex morphology, or move from $2D$ to $3D$. The addition of any one of these would increase computational time, in some cases exponentially. When modeling phenomena computationally, one lives by two questions: \textit{how much accuracy is required for validation?} and \textit{how can I live with that cost?}

%%%%%%%%%%%%%%%%%%%%%%%%%%%%%%%%%%%%%%%%%%%%%%%%
%
% QUANTITATIVE JELLYFISH CONVERGENCE: fluid
%
%%%%%%%%%%%%%%%%%%%%%%%%%%%%%%%%%%%%%%%%%%%%%%%%

\subsubsection{Error on Eulerian (fluid) grid}
\label{sec:jelly_eulerian_err}

Up to this point we have only discussed errors associated with the jellyfish's position, swimming speeds, and upward thrust force, without much mention of what is happening with the underlying fluid. That is, we've only explored errors associated with the Lagrangian (immersed structure). We have not discussed the fluid's velocity, pressure, nor vorticity. Unfortunately, discussing convergence of these quantities is not as cut and dry as those on the Lagrangian structure. For example, consider defining the absolute and relative errors analogously to Eqs.(\ref{eq:yForce_AbsErr})-(\ref{eq:yForce_RelErr}), for the fluid data, 

\begin{align}
    \label{eq:eulQ_AbsErr} err_{ABS_{\mbox{Q}}}^n &= \max_{ij} \left| \Upsilon Q_{1024_{ij}} - Q_{k_{ij}}  \right|,  \\
    \label{eq:eulQ_RelErr} err_{REL_{\mbox{Q}}}^n &= \max_{ij} \frac{ \left| \Upsilon Q_{1024_{ij}}  - Q_{k_{ij}} \right| }{ \left| \Upsilon Q_{1024_{ij}} \right| },
\end{align}

where $Q$ is a Eulerian quantity such as velocity, pressure, or vorticity, $n$ denotes the $n^{th}$ time-step, $k$ indicates a certain level of resolution below $1024\times1024$, and $\Upsilon$ is an interpolation operator for $1024\rightarrow k^{th}$ grid. Note that $Q$ is not spatially-averaged here.

\begin{figure}[H]
    %\centering
    \centering
    \includegraphics[width=0.99\textwidth]{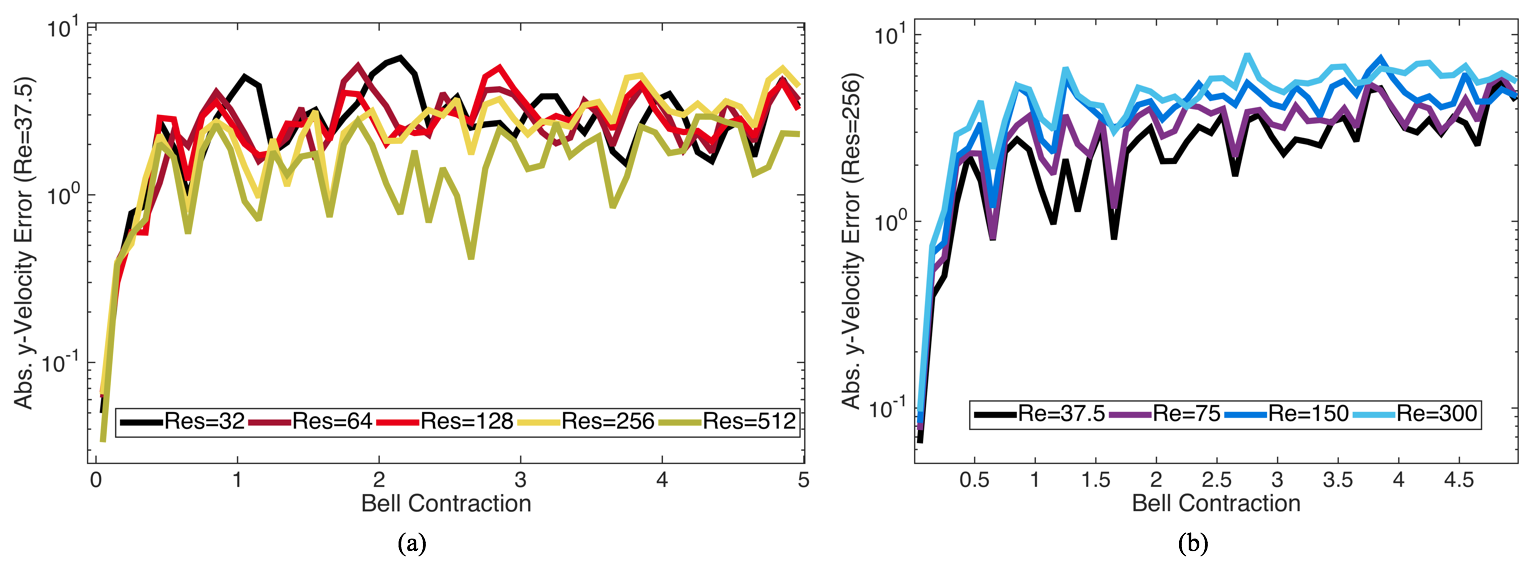}
    \caption{The absolute error of \textit{y-Velocity} across $5$ bell contractions when (a) $Re=37.5$ and grid resolution is varied and (b) when grid resolution is held at $256\times256$ and $Re$ is varied.}
    \label{fig:Jelly_AbsErr_Vs_Time_yVel}
\end{figure}

The absolute error of the fluid velocity in the y-direction (\textit{y-Velocity}) is shown in Figure \ref{fig:Jelly_AbsErr_Vs_Time_yVel} for $5$ bell contractions when either $Re$ is constant and grid resolution is varied (a) or vice-versa (b). In both cases errors start off small, increase, and appear to level-out. The higher grid resolution cases seem to lead to slightly less errors while higher $Re$ seems to tend toward slightly larger errors. However, are these showing any convergence? While higher resolution cases appear to have slightly less error, the absolute error appears to be time-dependent, where at some time-steps lower resolved cases actually have lower absolute error. How can this be? Perhaps it might be safe to say \textit{on average} higher resolved grids lead to less absolute error? 

\begin{figure}[H]
    %\centering
    \centering
    \includegraphics[width=0.99\textwidth]{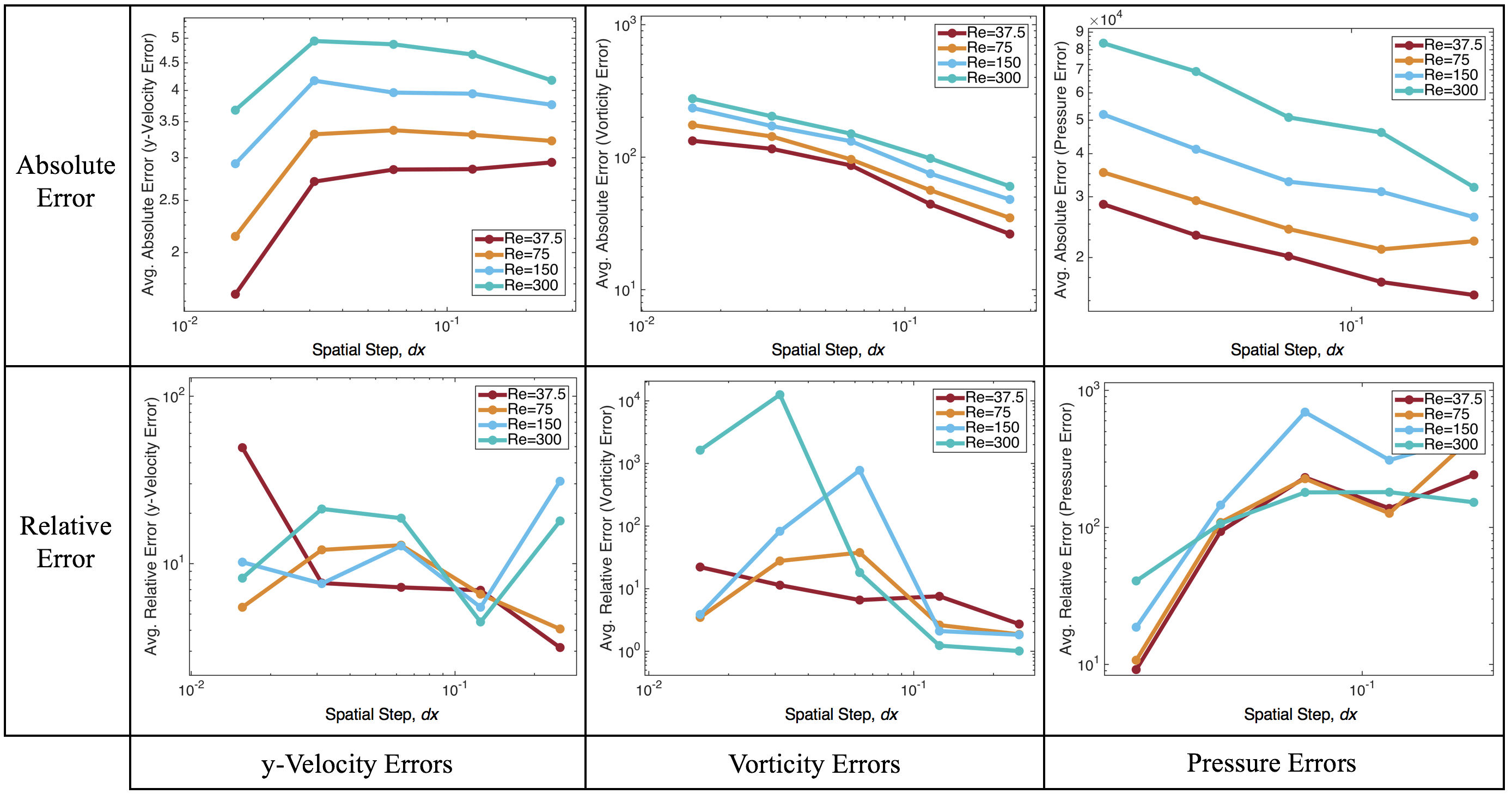}
    \caption{The temporally-averaged absolute and relative errors for \textit{y-velocity, vorticity}, and \textit{pressure} across $5$ bell contractions with varying $Re$.}
    \label{fig:Jelly_Eulerian_AbsRelErr}
\end{figure}

Figure \ref{fig:Jelly_Eulerian_AbsRelErr} gives the time-averaged absolute and relative errors for \textit{y-velocity, vorticity}, and \textit{pressure} across $5$ bell contractions with varying $Re$. In no case does a familiar looking convergence plot pop-up, illustrating higher grid-resolutions (smaller spatial steps) producing smaller errors. How can this be - are the simulations wrong? Let's investigate.

Recall Figure \ref{fig:Jelly_Conv_Sims_DistanceVsTime} in which shows how the same jellyfish swims different according to grid resolution for a particular $Re$ ($Re=150)$. Figures \ref{fig:Jelly_AbsErr_Vs_Time_yVel}-\ref{fig:Jelly_Eulerian_AbsRelErr} reflect this same phenomena, although it may be hidden at first. When the different resolution jellyfish are contracting, they do so in different regions of the domain as time moves forward, hence they are stirring up different areas within the fluid domain! Therefore the fluid dynamics between different resolution cases will be significantly different. This holds true even in the higher resolved cases where the jellyfish are near each other, but have slightly different kinematics. 

\begin{figure}[H]
    %\centering
    \centering
    \includegraphics[width=0.99\textwidth]{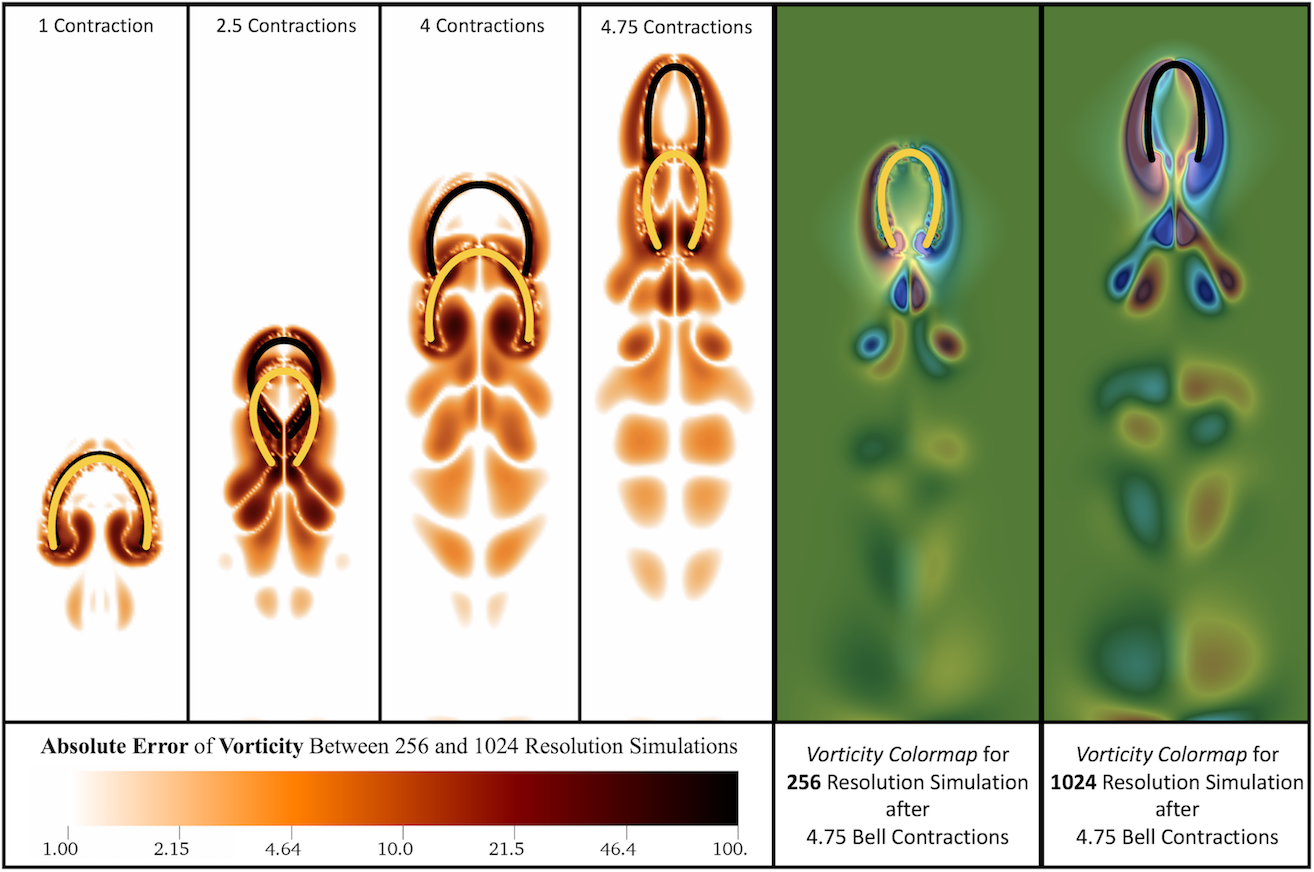}
    \caption{Colormaps of the absolute error of \textit{vorticity} between simulations for a $256\times256$ and $1024\times1024$ grid at different snapshots. The last error snapshot at $4.75$ bell contractions is compared to the actual vorticity colormaps for both grid resolution cases. These simulations are from a case of $Re=75$.}
    \label{fig:Jelly_Eulerian_Err_Vorticity}
\end{figure}

Figure \ref{fig:Jelly_Eulerian_Err_Vorticity} shows a spatial colormap of the absolute error of \textit{vorticity} over the course of a simulation when comparing grid resolutions of $256\times256$ and $1024\times1024$ for $Re=75$. It is clear that  vortex formation is significantly different between both cases due to different locations of the jellyfish, leading to a large absolute error in the vorticity. Moreover, from the last frame at $4.75$ bell contractions, it can be seen that some of the errors are magnified because of previously shed vorticies in the $1024$ case, which are oppositely spinning to those being created by the $256$ case at that moment. Similar behavior is seen in absolute errors for pressure and y-velocity, see Figures \ref{fig:Jelly_Eulerian_Err_yVel} and \ref{fig:Jelly_Eulerian_Err_Pressure} in Appendix \ref{App:EulerianError}, which present analogous absolute error data but for \textit{y-velocity} and \textit{pressure}, respectively. 

Simply put, since the jellyfish are swimming in different regions of the fluid domain, the fluid's behavior will be significantly different between differing resolution cases leading to magnification of errors! In fact, the errors tending towards a non-zero horizontal asymptote (Figure \ref{fig:Jelly_AbsErr_Vs_Time_yVel}) is an artifact of this. The highest absolute errors stem mostly from regions near each jellyfish bell. Once they are far enough apart the overall absolute errors are driven by the motion of each bell in a different region and when this occurs the absolute error steadies off. 

Since fluid-structure interaction systems errors traditionally display similar behavior at different grid resolutions, one way people can rectified this is to choose to only investigate the error up to a certain point in the simulation. For example Griffith et al. 2005 \cite{Griffith:2005} chose a time in which to compute errors when a deforming viscoelastic band completed one oscillation. Furthermore, rather than only look at averaged absolute or relative errors, one can define errors to be in terms of the $L^1$- and $L^2$-norm, or in general $L^p$-norm, e.g., 

\begin{equation}
    \label{eq:LP_Norm} \left|Q_{ij}\right|_p = \left( \sum_{ij} \left| \Upsilon Q_{1024_{ij}} - Q_{k_{ij}} \right|^p h^2 \right)^{1/p},
\end{equation}

where $p\geq 1$ (here $1$ or $2$), $Q$ is a scalar quantity, such as pressure, vorticity, or velocity in a particular direction, $k$ indicates a certain level of resolution below $1024\times1024$, $h$ is the spatial grid resolution corresponding to $k$, and $\Upsilon$ is an interpolation operator for $1024\rightarrow k^{th}$ grid. The $L^p$-norm is useful for calculating errors as no longer simply look for a maximal absolute or relative error at a certain time-step. Note that previously we had been using the $L^\infty$-norm to take the  maximal value of either absolute or relative error. 

\begin{figure}[H]
    %\centering
    \centering
    \includegraphics[width=0.99\textwidth]{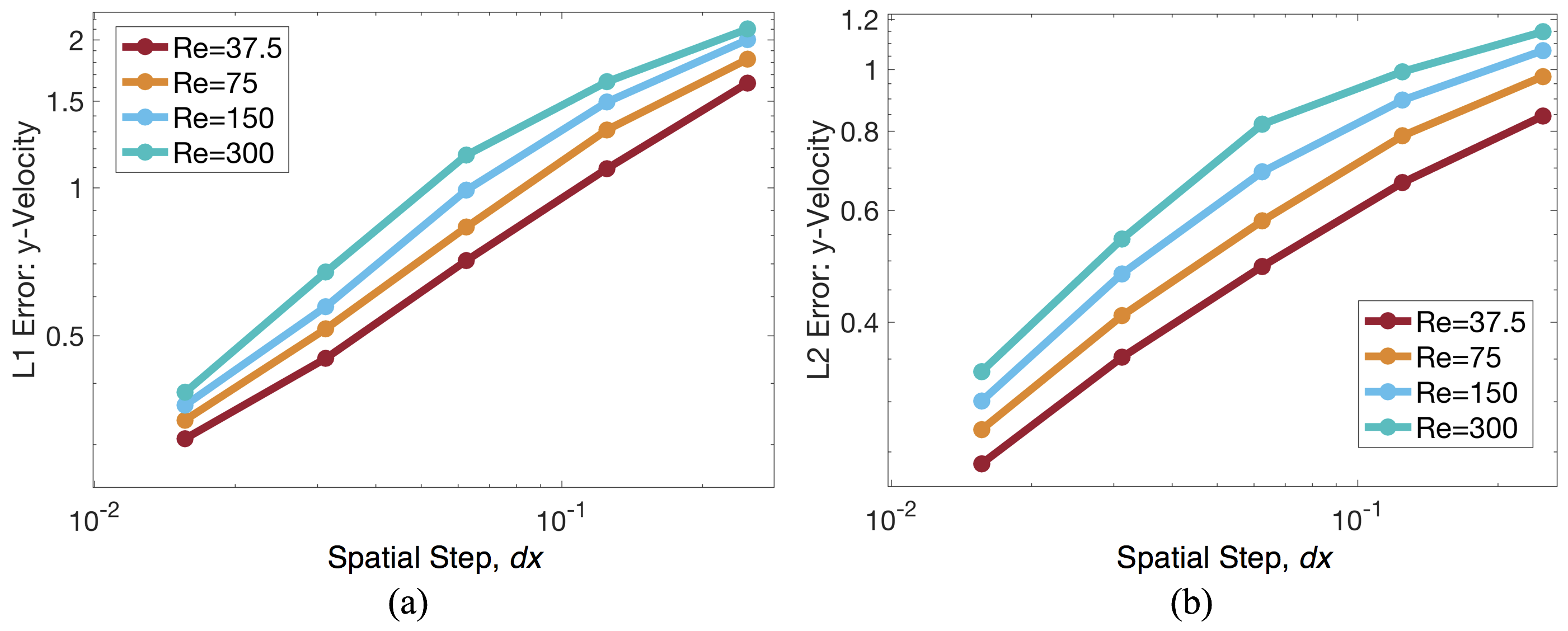}
    \caption{Depicting the (a) $L^1$- and (b) $L^2$-norms of the \textit{y-velocity} error after one bell contraction for various $Re$.}
    \label{fig:L1_L2_Error_differRe}
\end{figure}

Using Eq.(\ref{eq:LP_Norm}), we computed the $L^1$- and $L^2$- norms of the error for \textit{y-velocity} after the first bell contraction for various $Re$, as presented in Figure \ref{fig:L1_L2_Error_differRe}. The error decays as grid resolution increases (the spatial step-size decreases) under these norms; however, the errors do not appear small. Again, this is due to inconsistency with each jellyfish's location or kinematics after one bell contraction. The trend of lower $Re$ accompanied by lower error is still consistent with prior results. Moreover, the convergence rate (the slope of the line) is approximately consistent between all $Re$.

\begin{figure}[H]
    %\centering
    \centering
    \includegraphics[width=0.99\textwidth]{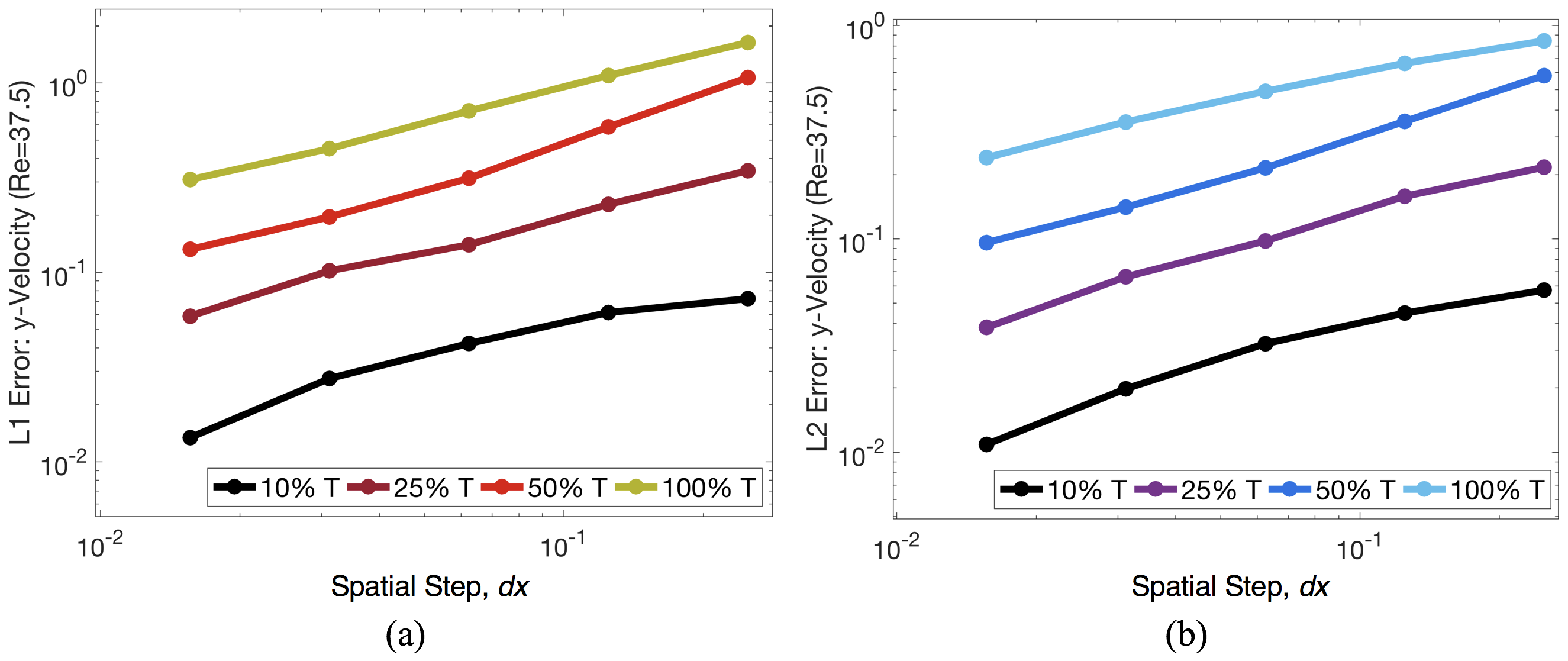}
    \caption{Depicting the (a) $L^1$- and (b) $L^2$-norms of the \textit{y-velocity} error during different percentages of one bell contraction period, $T$ for $Re=37.5$.}
    \label{fig:L1_L2_Error_differTime}
\end{figure}

Next we can investigate how the error magnitudes over one bell contraction for a particular $Re$. Figure \ref{fig:L1_L2_Error_differTime} shows how the error changes for different percentages of tbe first bell contraction, $T$. The errors increase further into the contraction cycle. For example at a grid resolution of $512\times512$ ($dx=0.015625$), at $10\% T$, both the $L^1-$ and $L^2$- errors are $\sim 0.01$, but by the end of a single bell contraction, they both are $\sim 0.5$. However, the convergence rate (the line's slope) remains approximately the same between all snapshots in time. 

To summarize, upon exploring the errors associated with the Eulerian (fluid) data in this fluid-structure interaction problems, we immediately observed that defining and interpreting errors was a non-trivial task, especially since there is no analytical solution in which to compare. In particular, we saw that (1) approaching the error calculations as we did previously for the same jellyfish model in Section \ref{sec:lagrangian_err} led to a question of whether we actually were observing convergence pertaining to the fluid data (Figure \ref{fig:Jelly_Eulerian_AbsRelErr}), (2) we had to rethink of how we thought about error for these simulations since there are objects moving around and pushing the fluid in different parts of the domain differently (Figures \ref{fig:Jelly_Eulerian_Err_Vorticity},\ref{fig:Jelly_Eulerian_Err_yVel} and \ref{fig:Jelly_Eulerian_Err_Pressure}), (3) we concluded we may only be able to compute errors in a semi-regular fashion during the first bell contraction (and by choosing an appropriate norm), otherwise the jellyfish would be too far away from each other in different resolution cases and it would be like comparing apples to oranges in the fluid domain (Figures \ref{fig:L1_L2_Error_differRe} and \ref{fig:L1_L2_Error_differTime}), and (4) simply striving for higher accuracy and boosting the resolution (decreasing the spatial step-size) could lead to simulations taking on the orders of weeks to run (Figure \ref{fig:Jelly_CompTime})! Moreover, it is unclear whether that boost (or time investment) in resolution would capture anymore biological/scientific relevance for the model being explored.

%MOVIES:
%-movie of jelly swimming at different resolutions
%-distance vs. time for swimming at different resolutions
%-emphasize same non-dimensional scaling coefficients for forces

%%%%%%%%%%%%%%%%%%%%%%%%%%%%%%%%%%%%%%%%%%%%%%%%%%%%%%%%%%%%%%%%%%%%%%%
%
% DISCUSSION AND CONCLUSIONS
%
%%%%%%%%%%%%%%%%%%%%%%%%%%%%%%%%%%%%%%%%%%%%%%%%%%%%%%%%%%%%%%%%%%%%%%%

%%%%%%%%%%%%%%%%%%%%%%%%%%%%%%%%%%%%%%%%%%%%%%%%
%
%
% DISCUSSION AND CONCLUSION
%
%
%%%%%%%%%%%%%%%%%%%%%%%%%%%%%%%%%%%%%%%%%%%%%%%%

\section{Discussion and Conclusion}
\label{sec:discussion}

Convergence plots are useful tools for both detailing numerical error for particular resolutions (number of iterations, time or spatial step-sizes) as well as for choosing appropriate resolutions to make practical use of an algorithm. Achieving minimal error tolerances is a quixotic endeavor if it take an unrealistic amount of computational resources and/or time to attain. In this paper we illustrated a few examples of convergence plots that stemmed from mathematical approximations to the Golden Ratio (Section \ref{sec:golden_ratio}), quadrature approximations using the Trapezoid Rule (Section \ref{sec:trap_rule}), basic ordinary differential equation time-stepping routines (Section \ref{intro} and Appendix \ref{app:euler}), and a contemporary research application involving jellyfish locomotion (Section \ref{jellyfish}). 

In Section \ref{sec:golden_ratio}, we introduced convergence plots as a way to determine how many terms were necessary in the Fibonacci Sequence to obtain certain levels of accuracy to numerical approximations to the Golden Ratio. In particular, $40$ terms were needed to achieve machine precision. Going beyond $40$ terms did not increase accuracy. From this example, the idea of maximal precision (\textit{double precision}) accuracy was introduced. Moreover, the computations performed in this application were trivial; they took virtually no time at all. 

Using the trapezoid rule from introductory calculus in Section \ref{sec:trap_rule}, convergence properties were shown to exhibit non-trivial behavior depending on the problem to which the numerical method is applied. Applying the composite trapezoid rule to a non-periodic function resulted in a slower convergence rate ($\sim 2-order)$) than when applied to an integrand that was periodic on its integration domain (\textit{exponential convergence}). In an example involving the former case, it takes more than $10^6$ sub-intervals to achieve machine precision where as in the latter, it took on the order of $10$. What a difference! This demonstrated that the same numerical method applied to different problems may result in significantly different accuracy or convergence properties!

Finally in Section \ref{jellyfish}, a contemporary fluid-structure interaction problem of a swimming jellyfish was used to illustrate use of convergence plots in a research and mathematical modeling setting. In particular, we saw that lower grid resolutions (and hence accuracy) resulted in decreased swimming performance, for both speed and distance traveled, for a range of different Reynolds Numbers. Moreover, investigating errors on the jellyfish itself (its position, swimming speed, and thrust force, e.g., Lagrangian data) was relatively straight-forward when comparing different grid resolution cases. However, exploring errors between fluid grids of different resolutions was non-trivial. This was mainly due to jellyfish contracting in different areas on the computational domain. %so naturally there would be large fluid dynamical differences between respective cases. 
%Furthermore, we posited that the amount of accuracy required for a mathematical model is non-trivial. 
For the case of jellyfish locomotion, increasing grid resolution could easily result in simulations requiring a week or longer of computing. On that note, it is unclear whether any higher grid resolution yields any additional benefits, i.e.,  captures more biological relevance. Instead, efforts may be better served in amending the model to capture more complex morphology, neuro-muscular paradigms, or other biological additions.

%%%%%%%This section feels repeated from above; maybe not needed?
Hopefully this has convinced you that there are significant benefits for performing convergence studies as well as some subtle nuances of numerical methods regarding both error and practicality. In particular, we demonstrated the following aspects of numerical methods and computational modeling:
\begin{enumerate}
    \item Additional resolution (e.g., more terms, a finer grid, etc.) does not achieve more accuracy if the method has already achieved machine precision accuracy (from Section \ref{sec:golden_ratio}).
    \item Applying the same numerical scheme to different problems can lead to drastically different convergence behavior (from Section \ref{sec:trap_rule}).
    \item Defining error itself may be non-trivial, e.g., for models with moving boundaries (from Section \ref{sec:quanitative_jellyfish}).
    \item Higher accuracy is not always practical or attainable if required computational time becomes unrealistic (from Sections \ref{intro} and \ref{jellyfish}).
    \item Higher accuracy may not always be the goal due to mathematical modeling assumptions (from Section \ref{jellyfish}), e.g., additional accuracy may not be important due to limitations of the model itself.
\end{enumerate}

Thus in practice for numerical simulation and mathematical modeling there is a trade-off between computational time, practicality, and desired accuracy. If a problem necessitates high accuracy but the computational time for a given numerical method is unfeasible, there are various actions one can take. One can either attempt a different numerical scheme or attempt to modify the existing method by implementing additional infrastructure to make the scheme faster, e.g., for differential equations some examples include possible adaptive time-stepping routines \cite{Atkinson:1989,Soderlind:2006}, adaptive mesh refinement (AMR) \cite{MJBerger84,MJBerger89}, parallelization \cite{Nievergelt:1964,Burrage:1993,Zhu:1994}, or even parallelization in time \cite{Emmett:2012,Gander:2013,Brandon:2015}. 
%Unfortunately both of those options are (almost certainly unpredictable) investments of one's time. 
This also tangentially evokes the practicality of novel numerical methods. When new methods are developed they are virtually always applied to \textit{toy} problems, meaning standard problems where solutions have been well-studied, categorized, and documented. Of course this is for the basis of comparison to existing methods and solutions; however, in many cases detailed convergence studies beyond those toy problems are non-existent, especially on contemporary research problems. Thus making it difficult to gauge how fruitful of one's time investment it may be to attempt to implement a new method. Fortunately, this is where interdisciplinary and integrative collaboration is key (and imperative)! 

The main purpose of this work was to immerse students in practical aspects of numerical analysis that border on contemporary research. While the notion of numerical error is stressed in all numerical analysis courses, it may be difficult for students to parse the subtleties that underlie error, such as those illustrated in this work. For this reason all codes, both simulation and analysis scripts are made available. As the scientific community becomes increasingly more dependent on mathematical modeling and numerical simulation, a thorough understanding of one's limitations due to speed, accuracy, and practicality is of the utmost importance.

\section*{Acknowledgments}
The authors would like to thank Charles Peskin for the development of immersed boundary method, Boyce Griffith for IBAMR, to which many of the input files structures of \textit{IB2d} are based, Alexander Hoover for supplying a IBAMR example of Jellyfish locomotion that was transcribed into \textit{IB2d} and made open source, and Laura Miller for continual comments, support, and conversation regarding mathematics and biology. They would also like to thank Christina Battista, Robert Booth, Namdi Brandon, Karen Clark, Jana Gevertz, Christina Hamlet, Christopher Jakuback, Andrea Lane, Jason Miles, Arvind Santhanakrishnan, Emily Slesinger, Christopher Strickland, and Lindsay Waldrop for comments on the design of the \textit{IB2d} software and suggestions for examples. NAB also wants to acknowledge his Spring 2018 Numerical Analysis class (Yaseen Ayuby, Shalini Basu, Gina Lee Celia, Rebecca Conn, Alexander Cretella, Robert Dunphy, Alyssa Farrell, Sarah Jennings, Edward Kennedy, Nicole Krysa, Aidan Lalley, Jessica Patterson, Brittany Reedman, Angelina Sepita, Nicole Smallze, Briana Vieira, and Ursula Widocki) for the original motivation for these project. This project was funded by the NSF OAC-1828163, TCNJ Support of Scholarly Activity (SOSA) Grant, the Department of Mathematics and Statistics, and the School of Science at TCNJ.
%\end{acknowledgements}

%%%%%%%%%%%%%%%%%%%%%%%%%%%%%%%%%%%%%%%%%%%%%%%%%%%%%%%%%%%%%%%%%%%%%%
%
% APPENDIX (if, appendix)
%
%%%%%%%%%%%%%%%%%%%%%%%%%%%%%%%%%%%%%%%%%%%%%%%%%%%%%%%%%%%%%%%%%%%%%%
%\clearpage
\appendix
%%%%%%%%%%%%%%%%%%%%%%%%%%%%%%%%%%%%%%%%%%%%%%%%
%
%
% APPENDIX: Euler Method 
%
%
%%%%%%%%%%%%%%%%%%%%%%%%%%%%%%%%%%%%%%%%%%%%%%%%

\section{Accuracy and Speed when using the Euler Method to solve ODEs}
\label{app:euler}

In Section \ref{intro} we introduced Figures \ref{fig:EulerComp} and \ref{fig:EulerTime} to illustrate how accuracy and speed are coupled, especially if a numerical scheme requires extensive computational time to run. 
%However, these plots were both introduced from a thousand foot view regarding numerical analysis and did not dive into the specific problem, numerical scheme, or application in which they arose. m
In this Appendix we explain the numerical method used to generate those plots.

These Figures came from solving an ordinary differential equation (ODE) with \textit{Euler's Method} applied tothe following ODE
\begin{equation}
    \label{eq:EulerODE} \frac{dy}{dt} = 2\pi\cos(2\pi t)
\end{equation}
with $y(0)=1$ over $0\leq t\leq 2$. 
%Now the right hand side of the above ODE is of the form in which we could simply use \textit{separable} methods to analytically solve it; however, we will use a computer, and in particular the Euler Method, to motivate the ideas of accuracy, speed, and convergence rates of numerical schemes.

The Euler Method can be derived in a few different ways, including stencils or interpolation \cite{Burden:2014}; however, Calculus students may favor the one in which we see the limit definition of a derivative come into fruition. Recall that the derivative of a function, say $y(t)$, at a point $t$ can be defined in the following manner,
\begin{equation}
    \label{eq:Derivative} \frac{dy}{dt}\Big|_{t} = \lim_{\Delta t\rightarrow 0} \frac{ y(t+\Delta t) - y(t) }{ \Delta t}.
\end{equation}

Since we are going to approximate the solution to ODEs using numerical methods, we note that a computer does not work in continua. This is precisely where (\ref{eq:Derivative}) will aid in solving (\ref{eq:EulerODE}). If we consider a fixed $\Delta t$ that is small \textit{enough}, we expect that the right hand side of (\ref{eq:Derivative}) could be at least a decent approximation of the derivative. This is the spirit of Euler's Method. Getting rid of the limit in the right hand side of (\ref{eq:Derivative}) and equating it with the right hand side of (\ref{eq:EulerODE}), we get the following:

\begin{equation}
    \label{eq:EulerMethod1} \frac{ y(t+\Delta t) - y(t) }{ \Delta t} = 2\pi\cos(2\pi t).  
\end{equation}

Solving for $y(t+\Delta t)$, we obtain a clear iterative form of Euler's Method,
\begin{equation}
    \label{eq:EulerMethod2} y(t+\Delta t) = y(t) + \Delta t \left( 2\pi\cos(2\pi t) \right).
\end{equation}

One can think of (\ref{eq:EulerMethod2}) as a recipe to step forward in time incrementally in \textit{time steps} of $\Delta t$. Since we know the initial value, $y(0)=1$, we have the first point in which we would step from. 

A natural question you may be wondering if \textit{how to choose $\Delta t$ so that it provides an accurate approximation of the derivative to even make this scheme meaningful?}. Setting aside any considerations of numerical stability (see \cite{Kincaid:2002,Burden:2014}), intuitively from seeing this derivation of the Euler Method, it is clear that the accuracy of this scheme is tied to the step size of $\Delta t$. This is precisely what the convergence plot given in Figure \ref{fig:EulerComp} from Section \ref{intro} illustrates. It is clear that as $\Delta t$ decreases, the accuracy increases; however, the price one has to pay for reducing $\Delta t$ is that the computational time required to run the simulation over the same window of $t$, $0\leq t\leq 2$, will significantly increase, see Figure \ref{fig:EulerTime}. 

%As mentioned in Section \ref{intro} there is a practical trade-off in the world of computing where one decides how much error their application can tolerate to understand what resolution (here $\Delta t$) to use. While striving towards minimal error (machine precision) is a noble pursuit, it could result in unfeasible (or unrealistic) computational time or computational resources.

The scripts used to perform this simulation and produce the Figures mentioned above is given in \textit{Supplemental/Eulers/}. %\textcolor{red}{https://youtu.be/SDVWr93Z9Yc?t=33s} \begin{verbatim} :) 
%\end{verbatim}

%%%%%%%%%%%%%%%%%%%%%%%%%%%%%%%%%%%%%%%%%%%%%%%%
%
%
% APPENDIX: Golden Ratio 
%
%
%%%%%%%%%%%%%%%%%%%%%%%%%%%%%%%%%%%%%%%%%%%%%%%%

\section{The Golden Ratio sequence is convergent and exhibits geometric convergence}
\label{app:fibonacci}

In this Appendix we prove that the sequence
\begin{equation*}
    \phi_n = \frac{F_{n}}{F_{n-1}}
\end{equation*}
converges to the Golden Ratio $\phi$. We begin with the following lemma.

\begin{lemma}\label{lem:append}
    For all $n\geq 1$ we have
    \begin{equation*}
        F_n^2-F_{n+1}F_{n-1} = (-1)^{n-1}.
    \end{equation*}
\end{lemma}

\begin{proof}
    It is clear for $n=1$. Inductively, assume
    \begin{equation*}
        F_n^2-F_{n+1}F_{n-1} = (-1)^{n-1}
    \end{equation*}
    for some $n\geq 1$. Then
    \begin{align}
        F_{n+1}^2 - F_{n+2}F_{n} &= F_{n+1}(F_n+F_{n-1})-(F_{n+1}+F_n)F_n \\
        &= F_{n+1}F_{n-1} +F_n^2 = -(-1)^{n-1} = (-1)^n.
    \end{align}*
        
\end{proof}

\begin{proof}[Proof of Convergence]
    Using Lemma \ref{lem:append} have
    
    \begin{align*}
        |\phi_n-\phi_{n-1}| &= \left| \frac{F_n}{F_{n-1}} - \frac{F_{n-1}}{F_{n-2}}\right| \\
        &= \left| \frac{F_{n}F_{n-1} -F_{n-1}^2}{F_n F_{n-2}}\right| \\
        &= \frac{1}{F_n F_{n-2}} \leq \frac{F_n-F_{n-2}}{F_n F_{n-2}} = \frac{1}{F_{n-2}}-\frac{1}{F_n}.
    \end{align*}
    So, by repeated application of the triangle inequality, with $n>m$:
    \begin{align} 
        |\phi_n-\phi_m| &\leq |\phi_n-\phi_{n-1}|+\cdots + |\phi_{m+1}-\phi_m| \label{cauchy_est1}\\
        &\leq \left(\frac{1}{F_{n-2}}-\frac{1}{F_n}\right) + \cdots  + \left(\frac{1}{F_{m-1}}-\frac{1}{F_{m+1}}\right) \nonumber\\
        &= \frac{1}{F_{n-2}}+\frac{1}{F_{m-1}}. \label{cauchy_est2}
    \end{align}
    Since $F_n\to \infty$ we have $|\phi_n-\phi_m|\to 0$ and so $\phi_n$ is Cauchy and thus convergent.
\end{proof}

{\bf Geometric convergence rate.}
Having justified the convergence of $\phi_n$, the analysis in Section \ref{sec:golden_ratio} shows $\phi_n\to \phi$. Take the limit $n\to \infty$ in \eqref{cauchy_est1}-\eqref{cauchy_est2} to obtain
\begin{equation}
    |\phi_m - \phi| \leq \frac{1}{F_{m-1}}.\label{geometric_rate}
\end{equation}
Since $\frac{F_{n+1}}{F_n} \to \phi$, then for any $\varepsilon>0$, we have (for large enough $n$), that $F_{n+1}> \phi F_n-\varepsilon$. As $\phi>1$, this implies that $F_{n}$ grows exponentially. 
Combining this with \eqref{geometric_rate}, it follows that the exponential growth rate of the Fibonacci sequence implies the geometric convergence rate of $\phi_n\to \phi$.

%%%%%%%%%%%%%%%%%%%%%%%%%%%%%%%%%%%%%%%%%%%%%%%%%%%%%%%%%%
%
%
% APPENDIX: SECANT METHOD CONVERGENCE
%
%
%%%%%%%%%%%%%%%%%%%%%%%%%%%%%%%%%%%%%%%%%%%%%%%%%%%%%%%%%%

\section{Secant Method Order of Convergence is the Golden Ratio}
\label{app:secant_convergence}

The \textit{Secant Method} is a popular method for computing roots of a non-linear equation, $f(x)=0$. It has a similar form to that of \textit{Newton's Method}, except it does not require an explicit evaluation of a derivative. Recall that the next root approximation for Newton's Method, $x_{n+1}$ is given by:
\begin{equation}
    \label{eq:NewtonsMetod} x_{n+1} = x_{n} - \frac{ f\left((x_n\right) }{ f'\left(x_n\right) }.
\end{equation}

Instead the Secant Method essentially approximates the derivative that found in the denominator of the fractional term in Eq.(\ref{eq:NewtonsMetod}), e.g., 

\begin{equation}
    \label{eq:derivApproxSecant} f'\left(x_n\right) \approx \frac{f\left(x_{n}\right)-f\left(x_{n-1}\right) }{x_{n}-x_{n-1}}.
\end{equation}

Note that while $x_n$ is the previous root approximation, $x_{n-1}$ is the second to last root approximation, and hence the Secant Method depends on the two previous approximations of the root. Therefore the next approximation of root of a nonlinear function by the Secant Method is given as

\begin{equation}
    \label{eq:SecantMethod} x_{n+1} = x_n - \frac{ f\left( x_n \right) }{\frac{f\left(x_{n}\right)-f\left(x_{n-1}\right) }{x_{n}-x_{n-1}}} = x_n - \frac{ f\left(x_n\right)(x_n-x_{n-1})}{ f\left(x_{n}\right)-f\left(x_{n-1}\right)}.
\end{equation}

\begin{figure}[H]
    %\centering
    \centering
    \includegraphics[width=0.50\textwidth]{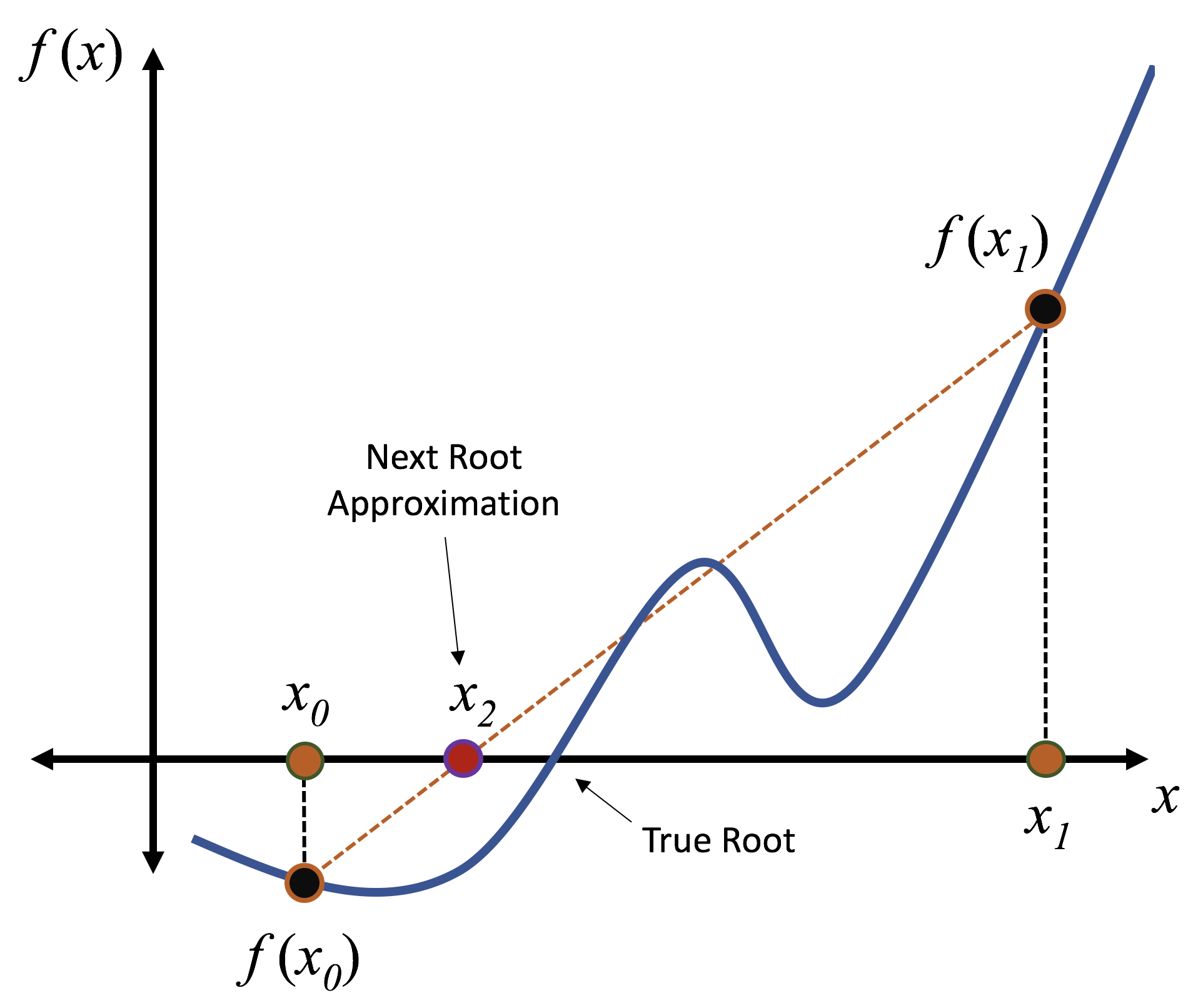}
    \caption{Graphical image depicting the derivation of the secant method for root finding, that is, computing the secant line between two successive guesses and choosing the next guess to be where it crossed the $x$-axis.}
    \label{fig:secant_derivation}
\end{figure}

However, there is no free lunch when using the Secant Method; the price we must pay for not computing an explicit derivative is that the Secant Method requires one additional initial guess and additional function evaluations. Furthermore, we will see that the order of convergence is also less than that of Newton's Method, which exhibits quadratic convergence for roots of multiplicity 1. Another derivation (and the one from whence it gets its name) is given graphically in Figure \ref{fig:secant_derivation}. 

The Secant Method can be derived by giving two initial guesses (previous root approximations), $x_0$ and $x_1$, finding the corresponding function values at those root approximations, $f(x_0)$ and $f(x_1)$, and then finding the secant line that connects those function values, e.g., 
$$y-f(x_1) = \left( \frac{f(x_1)-f(x_0)}{x_1-x_0} \right)(x-x_1).$$

Once we have the equation of that Secant Line, since we are searching for the root, we set $y=0$ and then solve for $x$, e.g.,
$$ x = x_2 = x_1 - \frac{ f(x_1) (x_1-x_0)}{f(x_1)-f(x_0)},$$
which is consistent with the form of Eq.(\ref{eq:SecantMethod}) above.

Next we wish to find the \textit{order of convergence} of the Secant Method. To do this we begin by assuming the root to which we're searching for is $x=\alpha$, e.g., $f(\alpha)=0$, and consider the error for the $(n-1)^{st}$, $n^{th}$, and $(n+1)^{st}$ root approximations, 
\begin{align}
    \nonumber e_{n-1} &= \alpha - x_{n-1}, \\
    \label{eq:SecantErrors} e_{n} &= \alpha - x_{n}, \\
    \nonumber e_{n+1} &= \alpha - x_{n+1}, 
\end{align}

where $\alpha$ is the true root. Next by substitution of Eqns.(\ref{eq:SecantErrors}) into (\ref{eq:SecantMethod}) we get 
\begin{equation}
    \label{eq:SecantErrorForm} e_{n+1} = \frac{ e_{n-1}f(x_n) - f(x_{n-1}) e_n }{f(x_{n}) - f(x_{n-1})}.
\end{equation}

The difficulty in the above equation is that the error at the next step, $e_{n+1}$, depends on both $e_{n}$ and $e_{n-1}$ as a linear combination of terms. To try to simplify this, we will invoke the Mean Value Theorem from Calculus that says in an interval between our current guess, $x_n$, and the true root, $\alpha$, there exists some number $\beta_n$ such that

$$f'(\beta_n) = \frac{ f(x_n) - f(\alpha) }{x_n-\alpha} = \frac{ f(x_n) }{x_n-\alpha} = \frac{ f(x_n) }{e_n},$$
since $f(\alpha)=0$ and $e_n=x_n-\alpha$ by Eq.(\ref{eq:SecantErrors}). Similarly, by Mean Value Theorem, there exists $\beta_{n-1}$ between $\alpha$ and $x_{n-1}$ such that 

$$f'(\beta_{n-1}) = \frac{ f(x_{n-1}) }{e_{n-1}},$$

and hence we can write

\begin{equation}
    \label{eq:secant_fs} f(x_n) = e_n f'(\beta_n) \ \ \ \mbox{ and } \ \ \ f(x_{n-1}) = e_{n-1} f'(\beta_{n-1}).
\end{equation}

The reason for invoking Mean Value Theorem will become obvious in a minute when we substitute Eq.(\ref{eq:secant_fs}) into (\ref{eq:SecantErrorForm}). However, while many times in mathematics we try to linearize problems, as nonlinearities are traditionally more difficult to study, in this particular problem we are actually trying extract a non-linear factor of $e_n e_{n-1}$ from both terms in Eq.(\ref{eq:SecantErrorForm}). Performing the aforementioned substitution we see

\begin{equation}
    \label{eq:secantErrorForm2} e_{n+1} = e_{n} e_{n-1} \frac{ f'(\beta_n)-f'(\beta_{n-1}) }{ f(x_n)-f(x_{n-1}) }.
\end{equation}

While the above does not look much different, it does have something special about it - namely the fractional term that is multiplying the errors is just some constant! Hence we see that the error at the $(n+1)^{st}$ step is proportional to the product of the errors at the previous two steps, e.g.,

\begin{equation}
    \label{eq:ErrorProportion} e_{n+1} \propto e_{n}e_{n-1}.
\end{equation}

Since we trying to determine the order of convergence, we can assume that the order is $p$, that is 

\begin{equation}
    \label{eq:orderOfConv} \lim_{n\rightarrow\infty} \frac{|x_{n}-\alpha|}{|x_{n-1}-\alpha|^p} = C \ \Longleftrightarrow \ e_{n} \propto (e_{n-1})^p,
\end{equation}
where $C$ is the asymptotic constant, which we will not concern ourselves with here. From Eq.(\ref{eq:orderOfConv}) we see that 

\begin{equation}
    \label{eq:SecantErrOrder1} e_{n+1} \propto (e_{n})^p
\end{equation}
by definition and by Eq.(\ref{eq:ErrorProportion}) previously that 

\begin{equation}
    \label{eq:SecantErrOrder2} e_{n+1} \propto e_{n}e_{n-1} \propto (e_{n-1})^p e_{n-1} = (e_{n-1})^{p+1}
\end{equation}

Hence combining Eqs.(\ref{eq:SecantErrOrder1}) and (\ref{eq:SecantErrOrder2}), we get that 

\begin{equation}
    \label{eq:SecantErrOrder3} (e_{n})^p \propto (e_{n-1})^{p+1} \ \Rightarrow \ e_{n} \propto (e_{n-1})^{\frac{p+1}{p}}.
\end{equation}
Using the relation from Eq.(\ref{eq:orderOfConv}) that $e_{n}\propto (e_{n-1})^p$ and substituting it into Eq.(\ref{eq:SecantErrOrder3}), we find that

$$(e_{n-1})^p \propto (e_{n-1})^{\frac{p+1}{p}},$$

and therefore gives us that $p = \frac{p+1}{p}$ which can be written in a more familiar quadratic form, e.g.,

\begin{equation}
    \label{eq:SecantLastStep} p^2-p-1=0.
\end{equation} 

This equation is the same equation we solved in Eq.(\ref{eq:GR_Quadratic}) from Section \ref{sec:golden_ratio}! Hence upon solving this quadratic we find that $$p=\frac{1\pm\sqrt{5}}{2}$$ and choosing the positive root, we obtain the Golden Ratio $\phi\approx 1.61803398874989485$. Therefore the Secant Method exhibits superlinear convergence, but not quite quadratic convergence!

%%%%%%%%%%%%%%%%%%%%%%%%%%%%%%%%%%%%%%%%%%%%%%%%%%%%%%%%%%
%
%
% APPENDIX: TRAP RULE EXPONENTIAL CONVERGENCE HEURISTIC
%
%
%%%%%%%%%%%%%%%%%%%%%%%%%%%%%%%%%%%%%%%%%%%%%%%%%%%%%%%%%%

\section{Trapezoid Rule: Heuristic Exponential Convergence with Periodic Functions}
\label{app:TrapRuleExpConv}

For simplicity of calculations, we assume that $f$ is a $2\pi$ periodic function with infinitely many bounded derivatives ($f\in C^\infty$). 
Periodic functions can be represented by the Fourier series: %which represents $f$ as a sum of trigonometric functions (as opposed to polynomials as in the Taylor series):
\begin{equation*}
    f(x) = \sum_{n=-\infty}^\infty c_n e^{inx},
\end{equation*}
where $e^{inx} = \cos(nx)+i \sin(nx)$, the famous Euler formula. The coefficients $c_n$ are given by 
\begin{equation*}
    c_n = \frac{1}{2\pi}  \int_0^{2\pi} f(x)e^{-inx}dx.
\end{equation*}

In particular, we have immediately that 
\begin{equation*}
    2\pi c_0 = I = \int_0^{2\pi} f(x) dx.
\end{equation*}
Applying the trapezoid rule to $f$ and using the Fourier series representation (and omitting details regarding absolute convergence which justifies switching sums):
\begin{align*}
    I_N &= \frac{2\pi}{2N}(f(0) + 2f(x_1)+\cdots + 2f(x_{N-1})+f(2\pi)) \\
    &= \frac{2\pi}{N} (f(0)+ f(x_1)+\cdots + f(x_{N-1})) \\
    &= \frac{2\pi}{N}\left( \sum_{n=-\infty}^{\infty} c_n e^{in\cdot 0} + \sum_{n=-\infty}^{\infty} c_n e^{in\cdot x_1} + \cdots + \sum_{n=-\infty}^{\infty} c_n e^{in\cdot x_{N-1}}\right) \\
    &= \frac{2\pi}{N} \sum_{m=0}^{N-1} \sum_{n=-\infty}^\infty c_n e^{in m \Delta x} \\
    &=\frac{2\pi}{N}  \sum_{n=-\infty}^\infty c_n\sum_{m=0}^{N-1} e^{in m \Delta x}.
\end{align*}
Since $\Delta x = \frac{2\pi}{N}$ we notice that
\begin{align*}
    c_n \sum_{m=0}^{N-1} e^{\frac{2\pi in m}{N}}  = 0
\end{align*}
whenever $n$ is not an integer multiple of $N$. Indeed, in that case $\{ e^{\frac{2\pi in m}{N}} \}_{m=0}^{N-1}$ generates the $N$th roots of unity on the complex plane (up to reordering) which sum to $0$ by symmetry.

On the other hand, if $n = kN$ for some $k\in \mathbb{Z}$, then
\begin{align*}
    c_{kN} \sum_{m=0}^{N-1} e^{2\pi i km} = N c_{kN}.
\end{align*}
Thus the trapezoid rule expression simplifies greatly:
\begin{equation*}
    I_N = 2\pi \sum_{k=-\infty}^\infty c_{kN}.
\end{equation*}
It remains to estimate
\begin{equation*}
    |I - I_N | = 2\pi \left| \sum_{k\neq 0} c_{kN} \right|.
\end{equation*}
To that end, integrate by parts $\alpha$ times (using periodicity of $f(x)$ and its derivatives) to obtain
\begin{align*}
    |c_{kN}| &= \left| \frac{1}{2\pi} \int_0^{2\pi} f(x) e^{-ikN x}dx \right| \\
    &=   \left|-\frac{1}{2\pi ikN} \int_0^{2\pi} f'(x) e^{-ikN x}dx + f(x)e^{-ikNx}\mid_{x=0}^{2\pi} \right| = \left| \frac{1}{2\pi kN} \int_0^{2\pi} f'(x) e^{-ikN x}dx \right|  \\
    &=\left|\frac{1}{2\pi k^2N^2} \int_0^{2\pi} f''(x) e^{-ikN x}dx\right|  \\
    &= \hspace{7mm} \vdots \\
    &=\left|\frac{1}{2\pi k^\alpha N^\alpha } \int_0^{2\pi} f^{(\alpha)}(x) e^{-ikN x}dx \right|.
\end{align*}
By assumption there is a constant $C$ so that for any $\alpha$, $\left|\int_0^{2\pi} f^{(\alpha)}(x) e^{-ikN x}dx\right|\leq C$, so
\begin{equation*}
    |I-I_N| \leq C \sum_{k\neq 0} \left|\frac{1}{k^\alpha N^\alpha} \right| = \frac{C}{N^\alpha} \sum_{k\neq 0} \left|\frac{1}{k^\alpha}\right|.
\end{equation*}
Since $\sum_{k\neq 0} k^{-\alpha}$ converges for any $\alpha>1$, we have $|I-I_N|\sim \frac{1}{N^\alpha}$. However since $\alpha$ was arbitrary, $I_N$ must converge to $I$ faster than any polynomial rate, hence suggesting geometric convergence. The ideas for this proof were motivated from \cite{Johnson:2010}.

%{\bf Faster convergence of of odd numbers of partitions.} We recall that the trapezoid rule appeared to converge faster for odd values of $N$ than for even values of $N$. To understand why this may occur, we rewrite
%\begin{align*}
%    \frac{1}{2\pi} |I - I_N | &=  \left|\sum_{k\neq 0} c_{kN} \right| = \left| \sum_{k=1}^\infty c_{kN}+c_{-kN} \right| \\
%    &= \left| \sum_{k=1}^\infty \int_0^{2\pi} f(x) (e^{ikNx}+e^{-ikNx}) dx \right| = 2\left| \sum_{k=1}^\infty \int_0^{2\pi} f(x) \cos(kNx)dx \right|. 
%\end{align*}
%Continuing heuristically, we see that if $N$ is even then all integrals contains terms of the form $\cos(2\ell x)$. In this case one expects ``resonance'' in the following sense: imagine a positive function $f$ which is supported predominantly around $\pi$. Then, the contribution to the infinite series from each term will necessarily be positive, since $\cos(2\ell \pi)=1$ for any $\ell$. On the other hand, if $N$ is odd then one expects ``additional cancellations:'' for the same function $f$ supported around $\pi$, the contribution to the infinite series from each term alternates between positive and negative values since $\cos(\ell \pi)=\pm 1$ depending on the parity of $\ell$. This additional cancellation results in smaller errors for odd values of $N$ than for even values. 

%%%%%%%%%%%%%%%%%%%%%%%%%%%%%%%%%%%%%%%%%%%%%%%%
%
%
% APPENDIX: Reynolds Number
%
%
%%%%%%%%%%%%%%%%%%%%%%%%%%%%%%%%%%%%%%%%%%%%%%%%

\section{The Reynolds Number}
\label{app:Re}

The Reynolds Number, $Re$, is a non-dimensional number in fluid dynamics that gives the ratio of inertial forces to viscous forces. One can think of viscous forces as those fluid forces that attempt to ``slow down" or ``impede fluid flow." In a nutshell, higher viscous forces arise from fluids with higher viscosities. High viscosity fluids include things like honey or corn syrup; fluids that are generally ``thicker" or "more sticky". The Reynolds Number does not only depend on viscosity and other physical characteristics of the fluid, e.g., its density), it also depends on characteristic length and velocity scales of the system being studied as well. 

The Reynolds Number, $Re$, is given quantitatively by the following expression
\begin{equation}
    \label{eq:Re} Re = \frac{\rho LV}{\mu},
\end{equation}
where $\rho$ and $\mu$ are the fluid's density and dynamical viscosity, respectively, while $L$ and $V$ are characteristic length and velocity scales for the system. Note that even if two systems may appear very different, if still have the same $Re$ they may display strikingly similar fluid behavior! For example, for us humans, swimming in water is generally no problem as long as we have prior swimming experience; however, for a bacteria trying to swim in water, it may feel like us trying to swim in peanut butter. In fact, the Re for us swimming in peanut butter is still approximately $5$ orders of magnitude greater than that of the swimming bacteria in water! E.g., 

$$Re_{\mbox{human}} = \frac{(\sim1280 kg/m^3)(\sim 1m)(\sim 1\ m/s)}{250 Pa\cdot s}\approx 7.7,$$

$$Re_{\mbox{bacteria}} = \frac{(1000 kg/s^2)(\sim 1m)(\sim 1\ m/s)}{0.001 Pa\cdot s}\approx 3\times10^-5,$$
of course assuming we are able to swim around $1\ m/s$ in peanut butter with densities and viscosities of $1280 kg/m^3$ and $250\ Pa\cdot s$, respectively \cite{Davis:2016}.

In the case of the Jellyfish, we define the characteristic velocity to be $fL$, the product of the speed of contraction and characteristic length scale, which is defined to be the width of the jellyfish at rest \cite{Hoover:2015}. The jellyfish model's non-dimensional parameters are found in Table \ref{table:params}.,

\begin{table}[H]
\begin{center}
\begin{tabular}{ |c|c|c| } 
 \hline
 Parameter Name & Symbol & Value   \\ \hline\hline
 Length Scale & $L$ & 1.0  \\ \hline
 Contraction Frequency & $f$ & 0.8  \\ \hline
 Fluid Density & $\rho$ & 1000 \\ \hline
 Fluid Dynamic Viscosity & $\mu$ & 6.66 \\ \hline
\end{tabular}
\end{center}
\caption{Non-dimensional parameters for the jellyfish model}
\label{table:params}
\end{table}

This gives $Re=150$, which is a biologically relevant $Re$ for jellyfish locomotion of certain jellyfish species \cite{Hershlag:2011,Hoover:2017}, such as \textit{Sarsia} \cite{Hoover:2015}.

For more information regarding Reynolds Numbers in biological applications please see \cite{Vogel:1996,Vogel:2013} or for $Re$ scaling research studies see \cite{Borazjani:2008,Hershlag:2011,Battista:2015b,Battista:2016a}.

%%%%%%%%%%%%%%%%%%%%%%%%%%%%%%%%%%%%%%%%%%%%%%%%%%%%%%%%%%%%%%%%%%%%%%
%
% APPENDIX: EXTRA JELLYFISH EULERIAN CONVERGENCE DATA
%
%%%%%%%%%%%%%%%%%%%%%%%%%%%%%%%%%%%%%%%%%%%%%%%%%%%%%%%%%%%%%%%%%%%%%%

\section{Additional Eulerian Absolute Error Data: \textit{y-Velocity} and \textit{Pressure}}
\label{App:EulerianError}

\begin{figure}[H]
    %\centering
    \centering
    \includegraphics[width=0.70\textwidth]{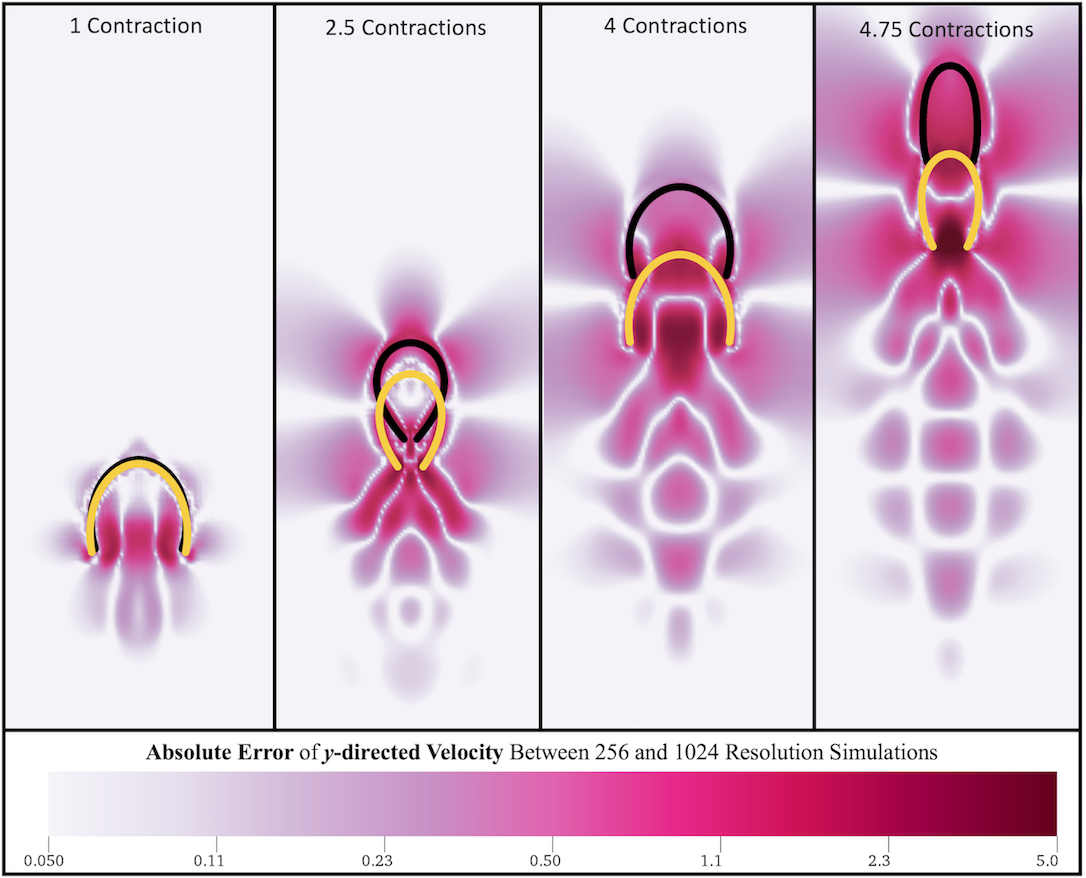}
    \caption{Colormaps of the absolute error of \textit{y-velocity} between simulations for a $256\times256$ and $1024\times1024$ grid at different snapshots for $Re=75$.}
    \label{fig:Jelly_Eulerian_Err_yVel}
\end{figure}

\begin{figure}[H]
    %\centering
    \centering
    \includegraphics[width=0.70\textwidth]{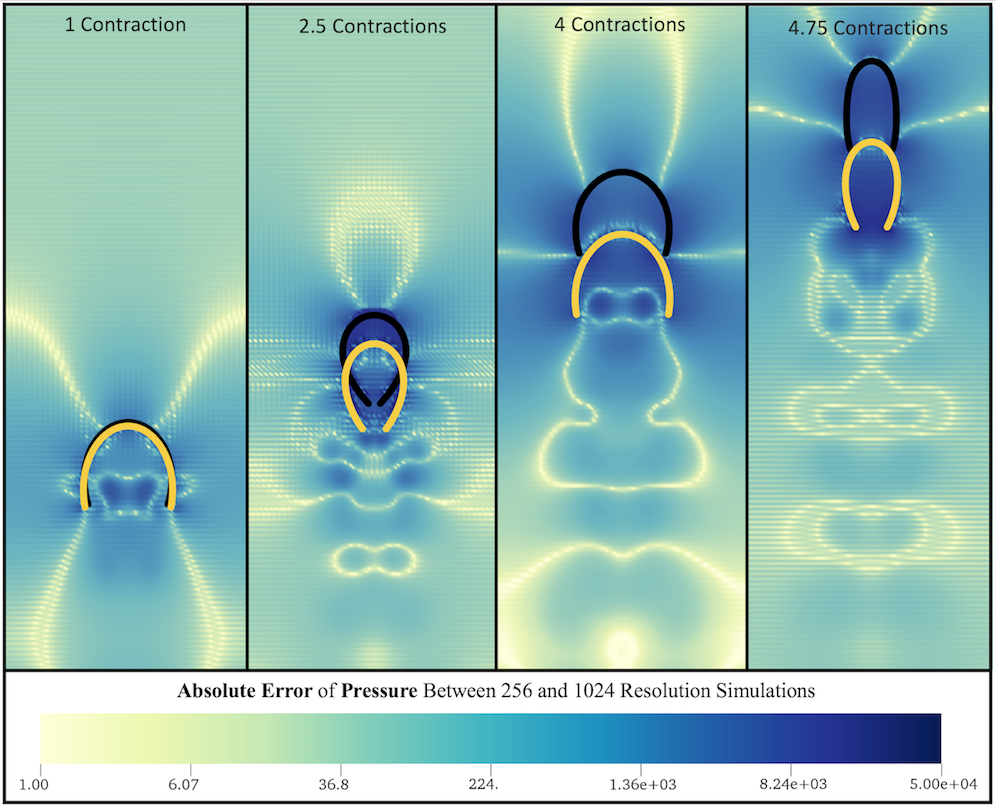}
    \caption{Colormaps of the absolute error of \textit{pressure} between simulations for a $256\times256$ and $1024\times1024$ grid at different snapshots for $Re=75$.}
    \label{fig:Jelly_Eulerian_Err_Pressure}
\end{figure}

Figures \ref{fig:Jelly_Eulerian_Err_yVel} and \ref{fig:Jelly_Eulerian_Err_Pressure} give a spatial depiction of the absolute error for \textit{y-velocity} and \textit{pressure}, respectively. 

%%%%%%%%%%%%%%%%%%%%%%%%%%%%%%%%%%%%%%%%%%%%%%%%%%%%%%%%%%%%%%%%%%%%%%
%
% APPENDIX: IMMERSED BOUNDARY METHOD
%
%%%%%%%%%%%%%%%%%%%%%%%%%%%%%%%%%%%%%%%%%%%%%%%%%%%%%%%%%%%%%%%%%%%%%%

\section{Details regarding \textit{IB2d} and the Immersed Boundary Method (IB)}
\label{IB:Appendix}

Here we will briefly introduce the fluid-structure interaction software used for computations, \textit{IB2d}, including the numerical method it hinges upon, the \textit{immersed boundary method} (IB).

%%%%%%%%%%%%%%%%%%%%%%%%%%%%%%%%%%%%%%%%%%%%%%%%%%%%%%%%%%%%%%%%%%%%%%
%
% APPENDIX: IB2d
%
%%%%%%%%%%%%%%%%%%%%%%%%%%%%%%%%%%%%%%%%%%%%%%%%%%%%%%%%%%%%%%%%%%%%%%

\subsection{\textit{IB2d}}
\label{app:IB2d_info}

Fluid dynamics, especially of the biological flavor, seems to be an ever encompassing subject. From the way organisms swim or fly to the way or our body helps us breath or even digest food, fluid dynamics, or more precisely, fluid-structure interactions (FSI) seem to be ever present. Traditionally this area of mathematical modeling is tightly wound with steep learning curves, which makes it challenging to teach effectively and give students meaningful first hand experiences. Fortunately, our open source software, \textit{IB2d}, was specifically designed to alleviate these challenges. It has two complete implementations in high-level programming environments, MATLAB and Python, which makes it accessible undergraduate and graduate students, and even scientists with limited programming experience.

As mentioned \textit{IB2d} was created for both educational and research purposes. It comes equipped with over $60$ built in examples that allow students to explore the world of fluid dynamics and fluid-structure interaction, from examples that illustrate standard fluid dynamics principles, such as flow around a cylinder for multiple Reynolds Numbers, or classical instabilities, like the Rayleigh-Taylor Instability, to examples that illustrate fully-coupled interactions of a fluid and flexible, deformable structures, including recent contemporary biologically motivated examples, such as aquatic locomotion or embryonic heart development. Some such examples are highlighted in \cite{Battista:2015,BattistaIB2d:2016,BattistaIB2d:2017}. This makes \textit{IB2d} suitable for either course projects or homework assignments, but also an ideal resource for a range of courses, ranging from mathematical modeling and mathematical biology courses to fluid mechanics to scientific computing/numerical analysis.

To aid in the mission above, multiple tutorial videos have been produced to help practictioners get familiar with the software:
\begin{itemize}
    %
    % TUTORIAL 1
    \item \textbf{Tutorial 1}: \url{https://youtu.be/PJyQA0vwbgU} \\
    \textit{An introduction to the immersed boundary method, fiber models, open source IB software, IB2d, and some FSI examples!}
    %
    % TUTORIAL 2
    \item \textbf{Tutorial 2}: \url{https://youtu.be/jSwCKq0v84s} \\
    \textit{A tour of what comes with the IB2d software, how to download it, what Example subfolders contain and what input files are necessary to run a simulation}
    %
    % TUTORIAL 3
    \item \textbf{Tutorial 3}: \url{https://youtu.be/I3TLpyEBXfE} \\
    \textit{The basics of constructing immersed boundary geometries, printing the appropriate input file formats, and going through these for the oscillating rubberband example from Tutorial 2}
    %
    % TUTORIAL 4
    \item \textbf{Tutorial 4}: \url{https://youtu.be/4D4ruXbeCiQ} \\
    \textit{The basics of visualizing data using open source visualization software called VisIt (by Lawrence Livermore National Labs). Using the oscillating rubberband from Tutorial 2 as an example to visualize the Lagrangian Points and  Eulerian Data (colormaps for scalar data and vector fields for fluid velocity vectors)}
\end{itemize}

Further details about \textit{IB2d}'s functionality can be found in \cite{Battista:2015,BattistaIB2d:2016,BattistaIB2d:2017}.

%%%%%%%%%%%%%%%%%%%%%%%%%%%%%%%%%%%%%%%%%%%%%%%%%%%%%%%%%%%%%%%%%%%%%%%
%
% IB APPENDIX: SUBSECTION: Governing Eqns of IB
%
%%%%%%%%%%%%%%%%%%%%%%%%%%%%%%%%%%%%%%%%%%%%%%%%%%%%%%%%%%%%%%%%%%%%%%%

\subsection{Governing Equations of IB}
\label{app:IB_eqns}

The conservation laws for momentum and mass that govern the motion of an incompressible, viscous fluid are written as the following set of coupled partial differential equations,

\begin{equation} %-----------------
   \rho\Big[\frac{\partial\textbf{U}}{\partial t}({\textbf{x} },t) +\textbf{U}({\textbf {x}},t)\cdot\nabla \textbf{U}({\textbf{x}} ,t)\Big]=  \nabla p({\textbf{x}},t) + \mu \Delta \textbf{U}({\textbf{x}},t) + \textbf{F}({\textbf{x}},t) \label{eq:NS1}
\end{equation}
\begin{equation}
      \nabla \cdot \textbf{U}({\bf x},t) = 0 \label{eq:NSDiv1}
\end{equation}
where $\textbf{U}({\textbf{x}},t)$, $p({\textbf{x}},t)$, and $\textbf{F}({\textbf{x}},t) $ are the fluid's velocity, pressure, and the force per unit area applied to the fluid by the immersed structure (e.g., the jellyfish), respectively. The quantities $\rho$ and $\mu$ are the fluid's density and dynamic viscosity, respectively. The system's independent variables are the time $t$ and the position ${\textbf{x}}$. The variables $\textbf{U}, p$, and $\textbf{F}$ are studied in an Eulerian framework on a fixed Cartesian mesh, $\textbf{x}$. We note that Eqs.(\ref{eq:NS1}) and (\ref{eq:NSDiv1}) are the conversation of momentum and mass equations for an incompressible, viscous fluid, respectively.

Deformations of the structure and the motion of the fluid are described by integral equations, known as the \textit{interaction} equations. This is one novelty of IB; the interaction equations translate all communication between the fluid (Eulerian) grid and immersed boundary (Lagrangian grid). They are written as the following integral equations containing delta function kernels,
\begin{align}
   {\bf F}({\bf x},t) &= \int {\bf f}(s,t)  \delta\left({\bf x} - {\bf X}(s,t)\right) dq \label{eq:force1} \\
   {\bf U}({\bf X}(s,t))  &= \int \textbf{U}({\bf x},t)  \delta\left({\bf x} - {\bf X}(s,t)\right) d{\bf x} \label{eq:force2}
\end{align}
where ${\bf f}(s,t)$ is the force per unit length applied by the boundary to the fluid as a function of Lagrangian position, $s$, and time, $t$. Note that $\delta({\bf x})$ is a three-dimensional delta function. The Cartesian coordinates of a material point labeled by the Lagrangian parameter, $s$, at time $t$, are given by ${\bf X}(s,t)$. The Lagrangian forcing term, ${\bf f}(s,t)$, describes the deformation forces of the immersed boundary at the material point labeled, $s$. Eq.(\ref{eq:force1}) spreads this force from the immersed boundary to the fluid through the external forcing term in Eq.(\ref{eq:NS1}). Eq.(\ref{eq:force2}) interpolates the velocity and moves the immersed boundary at the local fluid velocity, which enforces the no-slip condition. Each integral transformation uses a two-dimensional Dirac delta function kernel, $\delta$, to convert between Lagrangian variables and Eulerian variables.

As mentioned before, the use of delta functions as the kernel of integral equations Eqs.(\ref{eq:force1}-\ref{eq:force2}) is a main facet of IB. To approximate these integrals, discretized (and regularized) delta functions are used \cite{Peskin:2002}. We used one of the standard ones described in \cite{Peskin:2002}, e.g., $\delta_h(\mathbf{x})$, 
\begin{equation}
\label{delta_h} \delta_h(\mathbf{x}) = \frac{1}{h^3} \phi\left(\frac{x}{h}\right) \phi\left(\frac{y}{h}\right) \phi\left(\frac{z}{h}\right) ,
\end{equation}
where $\phi(r)$ is defined as
\begin{equation}
\label{delta_phi} \phi(r) = \left\{ \begin{array}{l} \frac{1}{8}(3-2|r|+\sqrt{1+4|r|-4r^2} ), \ \ \ 0\leq |r| < 1 \\    
\frac{1}{8}(5-2|r|+\sqrt{-7+12|r|-4r^2}), 1\leq|r|<2 \\
0 \hspace{2.1in} 2\leq |r|.\\
\end{array}\right.
\end{equation}

Deformation forces are calculated in a way that is specific to a particular application. For example, if the immersed boundary is allow to bend or stretch will determine what types of \textit{fiber models} (see \cite{BattistaIB2d:2016,BattistaIB2d:2017}) may be utilized to model the material properties of the structure. In this particular example of jellyfish locomotion, stiff springs are used to tether the Lagrangian mesh together and form the jellyfish bell, while springs with dynamically updating resting lengths are used to mimic the subumbrellar muscles that induce muscular contraction of the bell for the purpose of swimming. The jellyfish bell is able to retain its shape due to the inclusion of beams tethering Lagrangian points along the bell, which allow for bending but have a preferred curvature. In essence, springs allow for stretching and compressing and beams allow for bending of an immersed structure. Their deformation forces are described, respectively, as follows
\begin{align}
    \label{fiber_spring} \mathbf{F}_{spr} &= -k_{spr} \left( 1 - \frac{R_L}{\left|\left| \mathbf{X}_{S} - \mathbf{X}_M \right|\right| } \right) \cdot \left( \mathbf{X}_M - \mathbf{X}_S \right). \\
    \label{fiber_beam} \mathbf{F}_{beam} & =-k_{beam} \frac{\partial^4}{\partial s^4}\Big( \mathbf{X}(s,t) - \mathbf{X}_B(s,t) \Big),
\end{align}
where $k_{spr}$ and $k_{beam}$ are the spring stiffnesses and beam stiffnesses for springs and beams, respectively. The terms $X_{M}$ and $X_{S}$ in the spring forces represent the positions (in Cartesian coordinates) of the master and slave Lagrangian nodes at time, $t$. The parameter $R_L$ is a spring's resting length. The term $\mathbf{X}_B(s,t)$ in the bending force represents the preferred curvature of the configuration at time, $t$. 

We also include the use of \textit{target} points. Target points are used to either hold the geometry nearly rigid or prescribe the motion of an immersed structure. Here we use target points to create a wall that disrupts the flow in the top of the computational domain to reduce artifacts of periodic boundaries. Target points are modeled using a penalty force formulation and are written as the following,
\begin{equation}
{\bf f}(s,t) = k_{targ} \left(\textbf{Y}(s,t) - {\bf X}(s,t)\right),
\label{eq:force3}
\end{equation}
where $k_{targ}$ is a stiffness coefficient and $\textbf{Y}(s,t)$ is the prescribed position of the target boundary. Note that $\textbf{Y}(s,t)$ can be a function of both the Lagrangian parameter, $s$, and time, $t$; however, here they are static. For this model $k_{targ}$ was chosen to be large so that it would the top line of flow disruptors nearly rigid.

%

%%%%%%%%%%%%%%%%%%%%%%%%%%%%%%%%%%%%%%%%%%%%%%%%%%%%%%%%%%%%%%%%%%%%%%%
%
% IB APPENDIX: SUBSECTION: Numerical Algorithm
%
%%%%%%%%%%%%%%%%%%%%%%%%%%%%%%%%%%%%%%%%%%%%%%%%%%%%%%%%%%%%%%%%%%%%%%%

\subsubsection{Numerical Algorithm}
\label{app:IB_Numerical_Algorithm}

As stated in the main text, we impose periodic and no slip boundary conditions on a rectangular domain; however, through the use of target points, we reduce artifacts in the vertical velocity from periodicity. To solve Equations (\ref{eq:NS1}), (\ref{eq:NSDiv1}),(\ref{eq:force1}) and (\ref{eq:force2}) we update the velocity, pressure, position of the boundary, and force acting on the boundary at time $n+1$ using data from the previous time-step, $n$. The IB does this in the following steps \cite{Peskin:2002,BattistaIB2d:2016}:

\textbf{Step 1:} Computes the force density, ${\bf{F}}^{n}$ on the immersed boundary, from the current boundary configuration, ${\bf{X}}^{n}$.\\
\indent\textbf{Step 2:} Uses Eq.(\ref{eq:force1}) to spread these boundary deformation forces from the Lagrangian boundary mesh to the Eulerian (fluid) mesh. \\
\indent\textbf{Step 3:} Solves the Navier-Stokes equations, Equations (\ref{eq:NS1}) and (\ref{eq:NSDiv1}), on the Eulerian grid, which updates the fluid velocity, e.g., computing ${\bf{u}}^{n+1}$ and $p^{n+1}$ from ${\bf{u}}^{n}$, $p^{n}$, and ${\bf{f}}^{n}$.\\
\indent\textbf{Step 4:} Updates the Lagrangian positions, ${\bf{X}}^{n+1}$,  using the local fluid velocities, ${\bf{U}}^{n+1}$, computed from ${\bf{u}}^{n+1}$ and Equation (\ref{eq:force2}) via interpolation.

\bibliographystyle{siamplain}
\bibliography{heart}
\end{document}

\maketitle

%%%%%%%%%%%%%%%%%%%%%%%%%%%%%%%%%%%%%%%%%%%%%%%%%%%%%%%%%%%%%%%%%%%%%%
%
% SUPPLEMENT 1: THE MAIN SUPPLEMENTARY MATERIALS
%
%%%%%%%%%%%%%%%%%%%%%%%%%%%%%%%%%%%%%%%%%%%%%%%%%%%%%%%%%%%%%%%%%%%%%%

\section{Supplement 1:} The main supplementary file contains movies and codes pertaining to all the simulations detailed in this paper. It encompasses the following: $ $\\
\label{supp:1}

\begin{enumerate}
    \item \textbf{Convergence$\_$Classroom$\_$Supplement.pptx/.pdf}: presentations which may be used in class; slides that tell the story of the paper. Note that the $.pptx$ file has embedded movies in $.mp4$ format.
    \item \textbf{Eulers/}: directory containing the $.m$ file to run the ODE simulations with Euler's Method and produce the plots shown in the manuscript (Section 1 and Appendix A).
    \item \textbf{Golden$\_$Ratio/}: directory containing the $.m$ files to compute the Golden Ratio and produce the error plots shown in the manuscript (Section 2).
    \item \textbf{Trapezoid$\_$Rule/}: directory containing the $.m$ file to compute a definite integral using the Trapezoid Rule and produce the error plots shown in the manuscript (Section 3).
    \item \textbf{Jellyfish/}: directory containing the post-processed data, $.m$ error script files, simulation and Eulerian error movies, as well as simulation skeletons (source code and driving code) to produce all the simulation data presented in the manuscript (Section 4).
    \item \textbf{Projects/}: directory containing potential classroom activities or projects for students on Golden Ratio convergence (Section 2) or Trapezoid Rule convergence (Section 3).
\end{enumerate}

$ $ \\
We will now go into specifics about each sub-directory within this Supplemental Directory. $ $ \\

\begin{itemize}
    
    %
    % EULER'S
    %
    \item \textbf{Eulers/}: 
    \begin{itemize}
        \item \textit{Eulers.m}: script that performs Euler's Method of an ODE, computes associated error, and produces all the plots shown in the manuscript.
    \end{itemize}
    
    %
    % Golden Ratio
    %
    \item \textbf{Golden$\_$Ratio/}: 
    \begin{itemize}
        \item \textit{compute$\_$Golden$\_$Ratio$\_$Error.m}: script that Computes Golden Ratio Approximation for a Certain Number of Terms; produces plots given in the manuscript.
        \item \textit{find$\_$GoldenRatio$\_$Term$\_$Scaling.m}: script that Finds scaling relation for \# of terms in Fibonacci Sequence to approximate the Golden Ratio to a certain error tolerance; produces plots given in the manuscript
    \end{itemize}
    
    %
    % Trapezoid Rule
    %
    \item \textbf{Trapezoid$\_$Rule/}: 
    \begin{itemize}
        \item \textit{Trap$\_$Rule$\_$Error.m}: script that computes the error associated when applying Trapezoid Rule to  two different integral cases:
        \begin{enumerate}
            \item a periodic integrand on a periodic integration domain
            \item a non-periodic integrand;
        \end{enumerate}
        as well as produce all plots shown in the manuscript.
    \end{itemize}
    
    %
    % Jellyfish
    %
    \item \textbf{Jellyfish/}: 
    \begin{itemize}
        \item \textit{Simulation$\_$Skeletons/}: subfolder containing the source code (and driver files for thereof in \textit{IB2d$\_$Blackbox/} \cite{Battista:2015,BattistaIB2d:2016,BattistaIB2d:2017}) for all Jellyfish simulations presented in manuscript (e.g., for each Reynolds Number, $Re$, and Resolution considered).
        
        \item \textit{Data/}: the post-processed data from all the simulations corresponding to jellyfish position, forces acting upon, Eulerian data (velocity, vorticity, pressure), and their associated errors when compared to a high resolution case of $1024x1024$ grid resolution.
        \begin{enumerate}
            \item \textit{ReX$\_$Convergence$\_$Data/}: where $X=\{37pt5,75,150,300\}$ contains $.txt$ files where each corresponds to a particular simulation of a certain resolution. The data is stored in three columns, where the first is the simulation time, second is y-position of the top of the jellyfish bell, and third is instantaneous velocity of the top of the bell. NOTE: these files are read in by the script \textit{plot$\_$Position$\_$Velocity$\_$Data.m}.
            \item \textit{Forces$\_$Convergence$\_$Data}: folder containing $.txt$ files associated with the absolute and relative errors for each Reynolds Number, $Re$, studied. Each column in each $.txt$ file corresponds to a different grid resolution (resolution increases from left to right) and each row is a different time-step (first row is the first time-step stored, last row is the last time-step stored). NOTE: these files are read in by the script \textit{plot$\_$yForces$\_$Data.m}.
            \item \textit{Eulerian$\_$Convergence$\_$Data/}: folder contain $.txt$ files corresponding to the absolute and relative errors as well as $L^1$ and $L^2$ errors of the Eulerian Data (velocity, pressure, and vorticity) for each Reynolds Number, $Re$, considered. The column in each $.txt$ file corresponds to a different grid resolution, from lowest (left) to highest (right), while each row gives a different time-step. NOTE: these files are read in by the script \textit{plot$\_$Eulerian$\_$Data.m}.
            \item \textit{Eulerian$\_$Error$\_$VTK/}: folder containing a sub-folder, \textit{Error$\_$256} and a VisIt Session file \cite{HPV:VisIt} for visualizing the Eulerian error data. 
            \begin{itemize}
                \item \textit{Error$\_$256/}: contains $.vtk$ files that give the data for the absolute error of the Eulerian data (velocity, pressure, and vorticity) between the cases with $1024x1024$ and $256x256$ resolution for $Re=75$. It also contains two sub-folders that give the corresponding positions of the Lagrangian points in each case.
            \end{itemize}
        \end{enumerate}
            
        \item \textit{Error$\_$Scripts/}:
        \begin{enumerate}
            \item \textit{plot$\_$Position$\_$Velocity$\_$Data.m}: produces plots of the absolute and relative errors for the position and speed of the jellyfish; creates the plots given in the manuscript.
            \item \textit{plot$\_$yForces$\_$Data.m}: produces plots of the absolute and relative errors for the vertical forces acting on the jellyfish; creates the plots given in the manuscript.
            \item \textit{plot$\_$Eulerian$\_$Data.m}: produces plots of the absolute/relative and $L^1$/$L^2$ errors of the fluid (Eulerian) Data (velocity, pressure, and vorticity); creates the plots given in the manuscript.
        \end{enumerate}
        
        \item \textit{Movies/}:
        \begin{enumerate}
            \item \textit{Convergence$\_$Comparison/}: folder containing both $.mpg$ and $.mp4$ videos that illustrate jellyfish trying to swim at different Reynolds Numbers, $Re$, for varying grid resolutions. 
            \item \textit{Eulerian$\_$Error$\_$Movies/}: folder containing both $.mpg$ and $.mp4$ videos that illustrate the absolute error in the Eulerian data (velocity, pressure, vorticity) for jellyfish swimming at resolutions of $1024x1024$ and $256x256$ for $Re=75$. A VisIt Session file \cite{HPV:VisIt} is also included to reproduce the visualizations given here. It also contains a sub-folder (\textit{Swimmers$\_$1024-to-256/}) which has visualization movies (both $.mpg$ and $.mp4$) of each jellyfish swimming in the above examples and its associated Eulerian quantity (velocity, pressure, vorticity). 
        \end{enumerate}
        
    \end{itemize}

    %
    % Projects
    %
    \item \textbf{Projects/}: 
    \begin{itemize}
        \item \textit{Golden$\_$Ratio/}: folder contain a $.pdf$ (and corresponding LaTeX file) for a classroom activity or course project associated with Golden Ratio Convergence.
        \item \textit{Trapezoid$\_$Rule/}: folder contain a $.pdf$ (and corresponding LaTeX file) for a classroom activity or course project associated with Trapezoid Rule Convergence. 
    \end{itemize}

\end{itemize}

%%%%%%%%%%%%%%%%%%%%%%%%%%%%%%%%%%%%%%%%%%%%%%%%%%%%%%%%%%%%%%%%%%%%%%
%
% BIBLIOGRAPHY
%
%%%%%%%%%%%%%%%%%%%%%%%%%%%%%%%%%%%%%%%%%%%%%%%%%%%%%%%%%%%%%%%%%%%%%%

\bibliographystyle{siamplain}
\bibliography{heart}